\documentclass[11pt]{article}
\usepackage[nohead,margin=1.0in]{geometry}
\usepackage{amssymb, amsmath, amsthm}
\usepackage{color}
\usepackage{graphicx}
\usepackage{dsfont}
\usepackage{epstopdf}
\graphicspath{{./pics/}}
\usepackage{booktabs,siunitx}
\usepackage{threeparttable}
\usepackage{verbatim}
\usepackage{bm}
\usepackage[ocgcolorlinks, linkcolor=blue]{hyperref}

\newcommand{\bel}{\begin{equation} \label}
\newcommand{\ee}{\end{equation}}
\def\beq{\begin{equation}}
\def\eeq{\end{equation}}
\newcommand{\bea}{\begin{eqnarray}}
\newcommand{\eea}{\end{eqnarray}}
\newcommand{\beas}{\begin{eqnarray*}}
\newcommand{\eeas}{\end{eqnarray*}}

\newcommand{\R}{\mathbb{R}}
\newcommand{\C}{\mathbb{C}}

\newtheorem{theorem}{Theorem}[section]
\newtheorem{lem}[theorem]{Lemma}

\newtheorem{proposition}[theorem]{Proposition}

\numberwithin{equation}{section}

\newtheorem{thmx}{Theorem}

\renewcommand{\d}{\,\mathrm{d}}
\allowdisplaybreaks
\providecommand{\abs}[1]{\left\lvert#1\right\rvert}

\def\phi {\varphi}

\title{Simultaneous Identification of  Coefficients and Source in a Subdiffusion Equation from  One Passive Measurement\thanks{The work of B. Jin is
supported by Hong Kong RGC General Research Fund (Projects 14306423 and 14306824), ANR / RGC Joint Research
Scheme (A-CUHK402/24) and a start-up fund from The Chinese University of Hong Kong. The work of Y. Kian is supported by the French National Research Agency ANR and Hong Kong RGC Joint Research Scheme for the project IdiAnoDiff (grant ANR-24-CE40-7039).}}

\author{Maolin Deng\thanks{Department of Mathematics, The Chinese University of Hong Kong, Shatin, N.T., Hong Kong (\texttt{mldeng@link.cuhk.edu.hk, b.jin@cuhk.edu.hk})} \and Ali Feizmohammadi\thanks{ Department of Mathematical and Computational Sciences, University of Toronto Mississauga, 3359 Mississauga Road Deerfield Hall, Mississauga, ON L5L 1C6, Canada (\texttt{ali.feizmohammadi@utoronto.ca})}\and  Bangti Jin\footnotemark[2]\and
Yavar Kian\thanks{Univ Rouen Normandie, CNRS, Laboratoire de Math\'ematiques Rapha\"el Salem, UMR 6085, F-76000 Rouen, France (\texttt{yavar.kian@univ-rouen.fr})}}

\date{}

\begin{document}

\maketitle
\begin{abstract}
This article is devoted to the detection of parameters in anomalous diffusion from a single passive measurement. More precisely, we consider the simultaneous identification of coefficients as well as a time-dependent source term appearing in a time-fractional diffusion equation from a single boundary or internal passive measurement. We obtain several uniqueness results in dimension one as well as a multidimensional extension under some symmetry assumptions. Our analysis relies on spectral representation of solutions, complex and harmonic analysis combined with some known inverse spectral results for Sturm-Liouville operators. The theoretical results are complemented by a corresponding reconstruction algorithm and numerical simulations.  
\end{abstract}

\section{Introduction}

Diffusion is one of the most common transport mechanisms found in the nature, and the study of diffusion
processes occupies a central role in different scientific communities. In many practical applications, diffusion processes deviate from normal diffusion in the sense that the mean
squared particle displacement does not show a linear growth with time. Such diffusion
processes are collectively known as anomalous diffusion. The list of successful applications is still rapidly
growing and includes transport in complex media such as thermal diffusion in media with fractal geometry \cite{76}, dispersion in heterogeneous
aquifer \cite{1}, ion dispersion in column experiments \cite{18}, and protein transport in membranes \cite{57}; see
survey \cite{72} for an extended list of practical applications. A mathematically rigorous framework for capturing such phenomena is provided by \emph{subdiffusion equations}, where the standard time derivative is replaced by a \emph{fractional Caputo derivative of order} \( \alpha \in (0,1) \). This fractional operator encodes the memory and hereditary properties intrinsic to anomalous transport.

In recent years, significant attention has been devoted to \emph{inverse coefficient problems} \cite{CJLZ,CJYT,MillerYamamoto:2013,LiImanuvilovYamamoto:2016,HLYZ,JZ} or \emph{inverse source problems} \cite{JK3,JLLY,KLY23,KSXY,JK,Y} for subdiffusion equations, motivated both by their theoretical richness and their broad applicability to different problems, including
the identification of a source of diffusion
of pollution under the ground \cite{NSY}. Most existing studies (see the monograph \cite{KaltenbacherRundell:2023book} for details), however, operate under the assumption that either the medium is known or that one can \emph{actively interrogate} the system by prescribing boundary excitations, applying external forcing, or acquiring measurements over large regions of space-time. However, in many realistic settings, particularly in biomedical or remote sensing applications, such interventions are impractical or prohibited. In these scenarios, the observer must rely solely on \emph{passive measurements}, that is, recordings of the natural evolution of the system in response to unknown or uncontrolled input.

In this work, we address the inverse problem of determining a source of anomalous diffusion, described by one-dimensional subdiffusion, with unknown internal properties of the medium (e.g. unknown density of the medium or velocity of the moving quantities) from a \emph{single passive observation at a fixed spatial point away from the source over an infinite time sequence}. The source is determined simultaneously with the unknown properties of the medium associated with a coefficient of the subdiffusion equation. The problem can arise in determining a source of groundwater pollution diffusion \cite{Brusseau:1994,Sun:2024water} together with its velocity from a single passive measurement of the concentration of pollution away from the source. Our goal is to explore to what extent the  memory structure of subdiffusive dynamics --- manifested in the nonlocality of the fractional derivative --- permits the unique reconstruction of spatial information from such minimal data. This question draws conceptual inspiration from recent advances in the inverse analysis of classical evolution equations with passive measurements \cite{Feizmohammadi2025}, as well as from recent works devoted to the identification of an anomalous diffusion law in an unknown medium \cite{HongJinKian2024,JK1}.

The subdiffusive regime introduces profound analytical challenges. In contrast to the exponential decay behavior of parabolic equations, subdiffusion equations feature algebraically decaying Mittag-Leffler kernels, leading to a fundamentally different spectral representation of its solutions. Our analysis exploits these features to establish uniqueness results under minimal geometric and regularity assumptions, without relying on any active probing of the system.

To the best of our knowledge, this work provides the first theoretical investigation of passive inverse problems for time-fractional diffusion models. It contributes to the foundational understanding of how memory-driven dynamics can paradoxically \emph{enhance identifiability} in scenarios where data is scarce or observation is limited.
\subsection{Problem formulation and main results}
Let  $\ell>0$, $T>0$, $b\in C^1([0,\ell])$, and $q \in C([0,\ell])$. Given $\sigma\in L^2(0,T)$, $f\in L^2(0,\ell)$ and $\alpha \in (0, 1)$,
consider the unique weak solution $u$ of the following initial boundary value problem
\begin{equation}\label{eq1}
\left\{\begin{aligned}
\partial_t^{\alpha}u -\partial_x^2 u+b(x)\partial_xu+q(x)u &= \sigma(t)f(x), \quad \mbox{in }(0,\ell)\times(0,T),\\
 u(\ell,t)= 0,\ u(0,t)&= 0, \quad t\in(0,T), \\
u(x,0)&=0, \quad x\in(0,\ell).
\end{aligned}\right.
\end{equation}
Here, for all $\alpha\in(0,1)$, $\partial_t^\alpha$ denotes  the Caputo fractional derivative   of order $\alpha\in(0,1)$ defined by
$$\partial_t^\alpha h=D_t^\alpha (h-h(\cdot,0)),\quad h\in  C([0,+\infty);L^2(\Omega)),$$
where the Riemann-Liouville integral operator $I^\alpha$ and the Riemann-Liouville fractional derivative  $D_t^\alpha$, of order $\alpha$, are defined respectively by
\[
I^\alpha h(\cdot,t):=\frac{1}{\Gamma(\alpha)}\int_0^t\frac{h(\cdot,\tau)}{(t-\tau)^{1-\alpha}}\,\d\tau\quad \mbox{and}\quad D_t^\alpha:=
\partial_t\circ I^{1-\alpha}.
\]

For any Banach space $X$, we denote by $H_1(0,T;X)$ the space of functions $u\in H^1(0,T;X)$
satisfying the condition $u(0)=0$. For all $\beta\in (0,1)$, we denote by $H_\beta(0,T_1;X)$ the space of interpolation of order $1-\beta$ between $H_1(0,T;X)$ and $L^2(0,T;X)$ \cite[p. 12]{KRY}, which is a subspace of $H^\beta(0,T;X)$. In view of \cite[Theorem 2.1]{KRY}, problem \eqref{eq1} admits a unique solution $u\in H_\alpha(0,T;L^2(0,\ell))\cap L^2(0,T;H^{2}(0,\ell))$. Moreover, if $\sigma$ fulfills the extra assumption that there exists $\delta\in(0,T)$ such that
\bel{t1a}\sigma(t)=0,\quad t\in(T-\delta,T),\ee
then, $u|_{[0,\ell]\times(T-\delta,T]}\in C((T-\delta,T];C^1([0,\ell]))$, see Proposition \ref{p2} for the proof. With these observations for the forward problem, we are ready to state our inverse problem.

\begin{itemize}
    \item [({\bf IP})]{{\em Given a point $x_0 \in [0,\ell]$ and an increasing sequence $(t_n)_{n\in\mathbb N}$ in $(T-\delta,T)$ that converges to $T$, can we simultaneously and uniquely determine the coefficients $b$ and $q$ as well as the source function  $\sigma(t) f(x)$  from the knowledge of $\partial_xu(x_0,t_n)$ for $n\in\mathbb N$.}}
    \end{itemize}

For our first two results, we will assume that our measurements are made at an end point, say at $x_0=\ell$. Our first result assumes that the function $\sigma$ in \eqref{eq1} is a priori known.
\begin{thmx}\label{t1_intro} For $j=1,2$, let  $b_j\in C^1([0,\ell])$, $q_j\in C([0,\ell])$, let $f_j\in L^2(0,\ell)$ and let $\sigma\in L^2(0,T)$. Assume that $\sigma$ and $f_1$ are not identical to zero and that $\sigma$ satisfies \eqref{t1a}.
Define $V_j=-\frac{1}{2}b_j'+\frac{1}{4}b_j^2+q_j$ and assume that there exists $r\in (0,\ell)$ such that
$$
\mathrm{supp}(f_1)\subset [0,r] \quad \text{and}\quad \mathrm{supp} (V_1-V_2) \subset [0,\ell / 2-r / 2).$$
Assume (i) $q_1=q_2$ and $b_1(\ell)=b_2(\ell)$ or (ii) $b_1=b_2$.
Let $u_j$, $j=1,2,$ be the solution of \eqref{eq1} with  $b=b_j$, $q=q_j$, $f=f_j$. Then, given an arbitrarily increasing sequence $(t_n)_{n\in\mathbb N}$ in $(T-\delta,T)$ that converges to $T$, the condition
\begin{equation}\label{t1c}
\partial_xu_1(\ell,t_n)=\partial_xu_2(\ell,t_n),\quad \forall\, n\in\mathbb N
\end{equation}
implies that $b_1=b_2$, $q_1=q_2$ and $f_1=f_2$.
\end{thmx}

If $\alpha$ is not rational, we can extend the result to the determination of a more general class of time-dependent source term with separated variables.
\begin{thmx}\label{t2_intro} For $j=1,2$, let $\sigma_j\in L^2(0,T)$ be not identical to zero and assume that \eqref{t1a} is fulfilled with $\sigma=\sigma_j$ and denote by $u_j$ the solution of \eqref{eq1}, with  $\sigma=\sigma_j$, $b=b_j$, $f=f_j$ and $q=q_j$. Assume that the hypotheses of Theorem~\ref{t1_intro} is fulfilled. Then, condition \eqref{t1c} implies that
\bel{TT}b_1(x)=b_2(x),\quad q_1(x)=q_2(x),\quad \sigma_1(t)f_1(x)=\sigma_2(t)f_2(x),\quad (x,t)\in(0,\ell)\times(0,T).\ee
\end{thmx}
To the best of our knowledge, we obtain in Theorems \ref{t1_intro} and \ref{t2_intro} the first results in the literature of inverse problems for subdiffusion equations that enable simultaneous recovery of time-dependent source and space-dependent coefficients from a single passive boundary measurement of the corresponding evolution PDE. In fact, the only prior work dealing with the challenging inverse problem for the recovery of general coefficients of 1-D evolution equations from passive measurements can be found in \cite{Feizmohammadi2025}. The results of
Theorems \ref{t1_intro} and \ref{t2_intro} differ from those of \cite{Feizmohammadi2025} by the following aspects:\\
\begin{enumerate}
    \item [(i)] The work \cite{Feizmohammadi2025} is primarily concerned with recovery of coefficients and initial data for heat and wave equations from a single passive measurement. In this work, we consider subdiffusion equations and also address recovery of source terms which present some additional challenges.
    \item[(ii)] In contrast to \cite{Feizmohammadi2025}, Theorems \ref{t1_intro} and \ref{t2_intro} do not impose any condition on the source $f_2$ (the condition on support is only imposed on $f_1$).
    \item[(iii)] In Theorems \ref{t1_intro} and \ref{t2_intro}, we exploit the memory effect of subdiffusion to recover a general class of time-dependent sources from \textit{a posteriori} analysis of our single measurement.
\end{enumerate}

We would now like to consider {\bf (IP)} in the case of an interior measurement. Specifically, consider $a_1,a_2\in(0,\ell)$, with $a_1<a_2$, $a_1< \frac{\ell}{2}$ and $a_1+a_2\geq\ell$. By replacing the passive boundary measurement of the preceding results with an internal passive measurement in a neighborhood of $a_2$, together with some a priori knowledge of the coefficients on $[a_1,a_2]$, we have the following result.

\begin{thmx}\label{tt1_intro} For $j=1,2$, let  $b_j\in C^1([0,\ell])$, $q_j\in C([0,\ell])$. Consider $f_j\in L^2(0,\ell)$, $\sigma_j\in L^2(0,T)$, $j=1,2$, be a non-uniformly vanishing function such that  \eqref{t1a} is fulfilled with $\sigma=\sigma_j$ and denote by $u_j$ the solution of \eqref{eq1}, with  $\sigma=\sigma_j$, $b=b_j$, $f=f_j$ and $q=q_j$. Assume that
$$b_1'(x)+\frac{1}{4}b_1^2(x)+q_1(x)=b_2'(x)+\frac{1}{4}b_2^2(x)+q_2(x),\quad x\in[a_1,a_2]$$
and that one of the following conditions is fulfilled: (i) $q_1=q_2$ and $b_1(a_2)=b_2(a_2)$ or (ii) $b_1=b_2$. Moreover, if $\alpha\in(0,1)\cap \mathbb Q$ then we assume that $\sigma_1=\sigma_2$. Also, assume that there exists $r\in (0,\min(\ell-2a_1,2a_2-\ell))$ such that supp$(f_1)\subset [0,r]$.
Then,  given any increasing sequence $(t_n)_{n\in\mathbb N}$ in $(T-\delta,T)$ that converges to $T$ and any $\epsilon\in (0,\min(a_2,\ell-a_2))$, the condition
\bel{tt1a}u_1(x,t_n)=u_2(x,t_n),\quad n\in\mathbb N,\ x\in(a_2-\epsilon,a_2+\epsilon)\ee
implies that \eqref{TT} is true.
\end{thmx}

An important feature of Theorem \ref{tt1_intro} lies in the fact that it requires only the knowledge of the coefficients in the interval $[a_1,a_2]$ that can be arbitrarily small. This contrasts sharply with the results of Theorems \ref{t1_intro} and \ref{t2_intro}, as well as those of \cite{Feizmohammadi2025}, where the coefficients need to be known on at least half of the interval. The reason is that, contrary to the latter results, where the end point measurement at $x=\ell$ only produces access to the spectrum of a 1-D Sturm-Liouville operator, the case of a passive internal measurement also contains some nonlocal information about the associated eigenfunctions.

Finally, in the presence of domain symmetry, we have an extension of our results to multidimentional settings. Specifically, consider the cylindrical domain $\Omega=(0,\ell)\times\omega$ of $\R^n$, $n\geq2$ with $\omega$ a bounded domain of $\R^{n-1}$ with a $C^2$ boundary. Fix also $g\in L^2(\Omega)$ and consider the initial boundary value problem
\begin{equation}\label{eq111}
\begin{cases}
\partial_t^{\alpha}v -\Delta v+b(x_1)\partial_{x_1}u+q(x_1)u =  \sigma(t) g(x), & (x,t)=(x_1,x',t)\in (0,\ell)\times\omega\times(0,T),\\
 v= 0, & \mbox{on } \partial\Omega\times(0,T), \\
v=0, & \mbox{in } \Omega\times \{0\}.
\end{cases}
\end{equation}
By \cite[Theorem 2.1]{KRY}, we can show that this problem has a unique solution $u\in H_\alpha(0,T;L^2(\Omega))\cap L^2(0,T;H^{2}(\Omega))$. Moreover, in a similar way to Proposition \ref{p2}, we can prove that $w|_{\Omega\times(T-\delta,T]}\in C((T-\delta,T];H^2(\Omega)))$ (see \cite[Proposition 2.2]{JK1} for more details). Then, we can state our result as follows.

\begin{thmx}\label{t4_intro} For $j=1,2$,  $b_j\in C^1([0,\ell])$, $q_j\in C([0,\ell])$. For $j=1,2,$, let $g_j\in L^2(\Omega)$ and $\sigma_j\in L^2(0,T)$ be not identical to zero and assume that there exists $\delta\in(0,T)$ satisfying \eqref{t1a} with $\sigma=\sigma_j$. In addition, if $\alpha\in(0,1)\cap \mathbb Q$ then we assume that $\sigma_1=\sigma_2$. Let $V_j=-\frac{1}{2}b_j'+\frac{1}{4}b_j^2+q_j$ and assume that there exists $r\in(0,\ell)$ such that supp$(f_1)\subset [0,r]\times\omega$ and \eqref{t1b} is fulfilled. Finally, assume that one of the following conditions is fulfilled:
(i) $q_1=q_2$ and $b_1(\ell)=b_2(\ell)$;
(ii) $b_1=b_2$.
Denote by $v_j$ the solution of \eqref{eq111}, with  $\sigma=\sigma_j$, $b=b_j$, $g=g_j$ and $q=q_j$. Then the condition \eqref{t4a}
implies that $q_1=q_2$, $b_1=b_2$ and
$$ \sigma_1(t)g_1(x)=\sigma_2(t)g_2(x) \quad (x,t)\in \Omega\times (0,T).$$
\end{thmx}

To the best of our knowledge, in Theorem \ref{t4_intro} we obtain the first multidimensional result about simultaneous recovery of time-dependent source and coefficients. While the coefficients of the equation are subjected to the strong condition that they depend only on one space direction, the source terms are only subjected to an assumption related to their supports.

\subsection{Organization}
The rest of the article is organized as follows. We recall and prove several properties about representation and analytic extension in time of solutions of \eqref{eq1} in Section \ref{sec:prelim}. In Section \ref{sec:proof} we present the proofs of Theorems  \ref{t1_intro}, \ref{t2_intro}, \ref{tt1_intro} and \ref{t4_intro}. In fact, we will prove a slightly strengthened version of Theorems  \ref{t1_intro}, \ref{t2_intro}, \ref{tt1_intro} and \ref{t4_intro} that do not impose the condition that one of the two coefficients $b$ or $q$ is a priori known. In such a case, it is known that there is a natural obstruction to the simultaneous identification of both coefficients as well as the unknown source term. We recover the unknowns up to this optimal obstruction; see Theorems~\ref{t1}, \ref{t2}, \ref{tt1} and \ref{t4} for the precise formulations. Section \ref{sec:numer} will be devoted to numerical simulations based on the Levenberg-Marquardt method for the reconstruction of unknown parameters from a single passive measurement associated with the theoretical results. Finally, in the appendix we recall some inverse spectral results for one dimensional Schr\"odinger operators that will play an important role in our analysis.

\section{Preliminaries}\label{sec:prelim}
In this section we consider some representation properties of solutions of
\eqref{eq1} and we discuss regularity properties of its solutions such as analyticity in time. We denote by  $\left\langle \cdot, \cdot\right\rangle$ the scalar product, i.e.,
$\left\langle f,g \right\rangle=\int_0^\ell f(x)g(x)\d x$ for $f,g\in L^2(0,\ell).$
 Let $V=-\frac{1}{2}b'+\frac{1}{4}b^2+q$ and let $A$ be the unbounded operator acting on $L^2(0,\ell)$ with its domain $D(A)=H^1_0(0,\ell)\cap H^2(0,\ell)$ defined by
$$Af=-f''+Vf,\quad f\in D(A).$$
Recall  that $A$ is a selfadjoint operator whose  spectrum
consists of a strictly increasing and unbounded sequence of simple eigenvalues $(\lambda_{n})_{n\geq1}$. In the Hilbert space $L^2(0,\ell)$, for each eigenvalue $\lambda_n$, we fix an associated eigenfunction $\phi_{n}$ of unit $L^2(0,\ell)$-norm for the operator $A$. Then $ (\phi_n)_{n\in\mathbb N}$ forms an orthonormal basis in $L^2(0,\ell)$. We fix also
$$n_0=\min\{n\in\mathbb N:\ \lambda_n>0\}.$$
We need the two-parameter Mittag-Leffler function $E_{\alpha,\beta}(z)$ defined by (see, e.g., \cite[Section 1.2.1]{P} or \cite[p. 59]{Jin:2021book})
\begin{equation*}
  E_{\alpha,\beta}(z) = \sum_{n=0}^\infty \frac{z^n}{\Gamma(n\alpha+\beta)},\quad z\in \mathbb{C}.
\end{equation*}
Below we denote by $\tilde{f}$ the map given by
\begin{equation}\label{f}\tilde{f}(x) =e^{-\frac{1}{2}\int_\ell^xb(s)\d s}f(x),\quad x\in(0,\ell).\end{equation}
Consider the operator $L$ with domain $D(L)=H^1_0(0,\ell)\cap H^2(0,\ell)$ defined by
$$Lf=-f''+bf'+qf,\quad f\in D(L)$$
and observe that $L$ is the elliptic operator appearing in \eqref{eq1}. Notice that $L$ is not self-adjoint and therefore we cannot consider the classical spectral representation of solutions of \eqref{eq1} (see e.g. \cite{Jin:2021book,JZ,KY1,KRY}). Nevertheless, by using a change of variable, we can state a representation formula for solutions of problem \eqref{eq1} in the orthonormal basis $(\phi_n)_{ n\in\mathbb N}$.
\begin{proposition}\label{p1} Let $\sigma\in L^2(0,T)$, $b\in C^1([0,\ell])$ and $q\in L^\infty(0,\ell)$. Then the unique solution $u$ of problem \eqref{eq1} is given by
\bel{p1a}u(x,t)=e^{\frac{1}{2}\int_\ell^xb(s)\d s}\left(\sum_{k=1}^\infty \left(\int_0^t(t-s)^{\alpha-1}E_{\alpha,\alpha}(-\lambda_{k}(t-s)^{\alpha})\sigma(s)\d s\right)\langle \tilde{f},\phi_k\rangle\phi_{k}(x)\right),\ee
for all $ (x,t)\in (0,\ell)\times(0,T)$.
In addition, for a.e. $t\in(0,T)$, we have
\bel{p1b}\partial_xu(\ell,t)=\sum_{k=1}^\infty \left(\int_0^t(t-s)^{\alpha-1}E_{\alpha,\alpha}(-\lambda_{k}(t-s)^{\alpha})\sigma(s)\d s\right)\langle \tilde{f},\phi_k\rangle\phi_{k}'(\ell),\ee
where the sequence
\begin{equation*}
\sum_{k=1}^\infty \left(\int_0^t(t-s)^{\alpha-1}E_{\alpha,\alpha}(-\lambda_{k}(t-s)^{\alpha})\sigma(s)\d s\right)\langle \tilde{f},\phi_k\rangle\phi_{k}'(\ell)
\end{equation*}
converges  in the sense of a function lying in $L^2(0,T)$.
\end{proposition}
\begin{proof}  Let $\tilde{u}\in H_\alpha(0,T;L^2(0,\ell))\cap L^2(0,T;H^2(0,\ell))$ be defined by
$$\tilde{u}(x,t) =e^{-\frac{1}{2}\int_\ell^xb(s)\d s}u(x,t),\quad (x,t)\in (0,\ell)\times(0,T).$$
Then $\tilde{u}$ is the unique solution of
\begin{equation}\label{eq11}
\begin{cases}
\partial_t^{\alpha}\tilde{u} -\partial_x^2 \tilde{u}+V(x)\tilde{u} = \sigma(t)\tilde{f}(x), & \mbox{in }(0,\ell)\times(0,T),\\
 \tilde{u}(\ell,t)= 0,\ \tilde{u}(0,t)= 0 & t\in(0,T), \\
\tilde{u}(x,0)=0, & x\in(0,\ell).
\end{cases}
\end{equation}
By following the argument of \cite[Proposition 2.1]{JZ} (see also \cite[Theorem 2.4]{KY1}), we can prove
$$\tilde{u}(x,t)=\sum_{k=1}^\infty \left(\int_0^t(t-s)^{\alpha-1}E_{\alpha,\alpha}(-\lambda_{k}(t-s)^{\alpha})\sigma(s)\d s\right)\langle \tilde{f},\phi_k\rangle\phi_{k}(x),\quad (x,t)\in (0,\ell)\times(0,T)$$
which clearly implies \eqref{p1a}.
Following \cite[Theorem 2.4]{KY1}, since $\tilde{f}\in L^2(0,\ell)$ and $\sigma\in L^2(0,T)$, the sequence
\bel{p1d}\sum_{k=n_0+1}^N \left(\int_0^t(t-s)^{\alpha-1}E_{\alpha,\alpha}(-\lambda_{k}(t-s)^{\alpha})\sigma(s)\d s\right)\langle \tilde{f},\phi_k\rangle\phi_{k},\quad t\in(0,T),\ N\geq n_0+1\ee
converges in the sense of $L^2(0,T;D(A))\hookrightarrow  L^2(0,T;H^2(0,\ell))$. By the Sobolev embedding theorem, $H^2(0,\ell)$ is continuously embedded into $C^1([0,\ell])$. Thus, we deduce that $L^2(0,T;D(A))$ is continuously embedded in $L^2(0,T;C^1([0,\ell]))$ and we obtain
\bel{p1ddc}\begin{aligned}&\partial_x\left(\sum_{k=n_0+1}^\infty \left(\int_0^t(t-s)^{\alpha-1}E_{\alpha,\alpha}(-\lambda_{k}(t-s)^{\alpha})\sigma(s)\d s\right)\langle \tilde{f},\phi_k\rangle\phi_{k}\right)|_{x=\ell}\\
&=\sum_{k=n_0+1}^\infty \left(\int_0^t(t-s)^{\alpha-1}E_{\alpha,\alpha}(-\lambda_{k}(t-s)^{\alpha})\sigma(s)\d s\right)\langle \tilde{f},\phi_k\rangle\phi_{k}'(\ell),\quad t\in(0,T),
\end{aligned}\ee
with the sequence
\begin{equation*}
    \sum_{k=n_0+1}^\infty \left(\int_0^t(t-s)^{\alpha-1}E_{\alpha,\alpha}(-\lambda_{k}(t-s)^{\alpha})\sigma(s)\d s\right)\langle \tilde{f},\phi_k\rangle\phi_{k}'(\ell)
\end{equation*}
converging  in the sense of a function lying in $L^2(0,T)$.
In addition, recalling that the map $E_{\alpha,\alpha}$ is holomorphic in $\mathbb C$, one can easily check that the map
$$\sum_{k=1}^{n_0} \left(\int_0^t(t-s)^{\alpha-1}E_{\alpha,\alpha}(-\lambda_{k}(t-s)^{\alpha})\sigma(s)\d s\right)\langle \tilde{f},\phi_k\rangle\phi_{k},\quad t\in(0,T)$$
is lying in $ L^2(0,T;H^2(0,\ell))\hookrightarrow  L^2(0,T;C^1([0,\ell]))$ and the map
$$\begin{aligned}&\partial_x\left(\sum_{k=1}^{n_0} \left(\int_0^t(t-s)^{\alpha-1}E_{\alpha,\alpha}(-\lambda_{k}(t-s)^{\alpha})\sigma(s)\d s\right)\langle \tilde{f},\phi_k\rangle\phi_{k}\right)|_{x=\ell}\\
&=\sum_{k=1}^{n_0} \left(\int_0^t(t-s)^{\alpha-1}E_{\alpha,\alpha}(-\lambda_{k}(t-s)^{\alpha})\sigma(s)\d s\right)\langle \tilde{f},\phi_k\rangle\phi_{k}'(\ell),\quad t\in(0,T)
\end{aligned}$$
is lying in $L^2(0,T)$. By
combining this with \eqref{p1ddc}, we deduce that
 the sequence \eqref{p1d} converges in the sense of $L^2(0,T;C^1([0,\ell]))$ and  get
$$\partial_x\tilde{u}(\ell,t)=\sum_{k=1}^\infty \left(\int_0^t(t-s)^{\alpha-1}E_{\alpha,\alpha}(-\lambda_{k}(t-s)^{\alpha})\sigma(s)\d s\right)\langle \tilde{f},\phi_k\rangle\phi_{k}'(\ell),\quad t\in(0,T),$$
with the sequence
\begin{equation*}
    \sum_{k=1}^\infty \left(\int_0^t(t-s)^{\alpha-1}E_{\alpha,\alpha}(-\lambda_{k}(t-s)^{\alpha})\sigma(s)\d s\right)\langle \tilde{f},\phi_k\rangle\phi_{k}'(\ell)
\end{equation*}
converging  in the sense of a function lying in $L^2(0,T)$. Finally, using the fact that
$$\partial_x\tilde{u}(\ell,t)=-\tfrac{1}{2}b(\ell)u(\ell,t)+\partial_xu(\ell,t)=\partial_xu(\ell,t),\quad t\in(0,T),$$
we deduce that the identity \eqref{p1b} holds. This completes the proof of the proposition.
\end{proof}

We also have the following result about the analytic extension of the solution of problem \eqref{eq1} under condition \eqref{t1a}.
\begin{proposition}\label{p2} Let the condition of Proposition \ref{p1} be fulfilled with $\sigma$ satisfying \eqref{t1a}. Then the restriction of the map $t\mapsto\partial_xu(0,t)$ to $(T-\delta,T)$ admits an analytic extension to $(T-\delta,+\infty)$ given by the map
\bel{p2a}h(t)=\sum_{k=1}^\infty \left(\int_0^{T-\delta}(t-s)^{\alpha-1}E_{\alpha,\alpha}(-\lambda_{k}(t-s)^{\alpha})\sigma(s)\d s\right)\langle \tilde{f},\phi_k\rangle\phi_{k}'(\ell),\quad t\in(T-\delta,+\infty).\ee
Moreover, fixing $\mathbb C_{n_0}=\{z\in\mathbb C: |z|\not\in \{|\lambda_n|^{\frac{1}{\alpha}}:\ n\in\mathbb N\}\}\setminus (-\infty,0]$,   the Laplace transform $\hat{h}(p)$ of $h$  admits an holomophic extension to $p\in \mathbb C_{n_0}$ given by
\bel{p2b}\hat{h}(p)=\hat{\sigma}(p)\sum_{k=1}^\infty\frac{\langle \tilde{f},\phi_k\rangle\phi_{k}'(\ell)}{p^\alpha+\lambda_k},\ee
where $\hat{\sigma}(p)$ is the Laplace transform of $\sigma$ extended by zero to $(0,+\infty)$ given explicitly by
$\hat{\sigma}(p)=\int_0^T\sigma(t)e^{-pt}\d t$ for $p\in\mathbb C.$

\end{proposition}
\begin{proof} In view of \eqref{p1b} and  \eqref{t1a}, we have
\begin{align*}
\partial_xu(0,t)&=\sum_{k=1}^\infty \left(\int_0^t(t-s)^{\alpha-1}E_{\alpha,\alpha}(-\lambda_{k}(t-s)^{\alpha})\sigma(s)\d s\right)\langle \tilde{f},\phi_k\rangle\phi_{k}'(\ell)\\
&=\sum_{k=1}^\infty \left(\int_0^{T-\delta}(t-s)^{\alpha-1}E_{\alpha,\alpha}(-\lambda_{k}(t-s)^{\alpha})\sigma(s)\d s\right)\langle \tilde{f},\phi_k\rangle\phi_{k}'(\ell)=h(t),\quad t\in(T-\delta,T).
\end{align*}
For any $\beta\in(0,\pi)$, $\mathcal C_\beta$ denotes the cones in the complex plane defined by $\mathcal C_\beta=\{-re^{i\theta}:\, r>0,\ \theta\in(-\beta,\beta)\}$ and $\mathcal C_{T-\delta,\beta}=\{T-\delta +re^{i\theta}:\, r>0,\ \theta\in(-\beta,\beta)\}$. From \cite[pp. 34--35]{P}, we know that there exists $\beta\in(0,\pi)$ such that
\[
\abs{E_{\alpha,\alpha}(z)}\leq\frac{C}{1+|z|},\quad z\in \mathcal C_\beta.\]
From this estimate, we deduce that there exists $\beta_\star \in(0,\pi)$ such that
\bel{pp2b} \abs{E_{\alpha,\alpha}(-\lambda_k(z-s)^\alpha)}\leq C\lambda_k^{-1}|z-T-\delta|^{-\alpha},\quad z\in \mathcal C_{T-\delta,\beta_\star},\quad k\in\mathbb N,\ s\in(0,T-\delta).\ee
Using these properties in a similar way to \cite[Proposition 2.1]{KSXY}, we can define the map
\begin{equation*}
F_1(z)=\sum_{k=n_0+1}^\infty \left(\int_0^{T-\delta}(z-s)^{\alpha-1}E_{\alpha,\alpha}(-\lambda_{k}(z-s)^{\alpha})\sigma(s)\d s\right)\langle \tilde{f},\phi_k\rangle\phi_{k}'(\ell),\quad z\in \mathcal C_{T-\delta,\beta_\star}
\end{equation*}
and proves that it is holomorphic with respect to $z\in \mathcal C_{T-\delta,\beta_\star}$. In addition, recalling that the map $E_{\alpha,\alpha}$ is holomorphic in $\mathbb C$, one can easily check  that the map
$$F_2(z)=\sum_{k=1}^{n_0} \left(\int_0^{T-\delta}(z-s)^{\alpha-1}E_{\alpha,\alpha}(-\lambda_{k}(z-s)^{\alpha})\sigma(s)\d s\right)\langle \tilde{f},\phi_k\rangle\phi_{k}'(\ell)$$
is analytic with respect to $z\in \mathcal C_{T-\delta,\beta_\star}$ and the same will be true for $F(z)=F_1(z)+F_2(z)$. Since
$h(t)=F(t)$ for  $t\in(T-\delta,+\infty)$,
we deduce that the restriction of the map $h$ to $(T-\delta,+\infty)$ is analytic. The last statement of the proposition can be deduced from  arguments similar to \cite[Proposition 2.2]{JK1} combined with the above argumentation.
\end{proof}

\section{Proofs of the main results}
\label{sec:proof}
In this section, we present the proofs of the main results.
\subsection{Proof of Theorem \ref{t1_intro}}\label{ssec:thm1}
We prove a slightly stronger version of Theorem~\ref{t1_intro} as follows.
\begin{theorem}\label{t1} For $j=1,2$, let  $b_j\in C^1([0,\ell])$, $q_j\in C([0,\ell])$. Consider $f_j\in L^2(0,\ell)$, $j=1,2$,  $\sigma\in L^2(0,T)$ and assume that $\sigma$ and $f_1$ are non-uniformly vanishing with $\sigma$ satisfying \eqref{t1a}.
Fix also $V_j=-\frac{1}{2}b_j'+\frac{1}{4}b_j^2+q_j$ and assume that there exists $r\in (0,\ell)$ such that $\mathrm{supp}(f_1)\subset [0,r]$ and
\begin{equation} \label{t1b}
\mathrm{supp} (V_1-V_2) \subset [0,\ell/2-r/2).
\end{equation}
Then, by denoting by $u_j$ the solution of \eqref{eq1}, with  $b=b_j$, $q=q_j$ and $f=f_j$, we have $u_j|_{[0,\ell]\times(T-\delta,T]}\in C((T-\delta,T];C^1([0,\ell]))$ and, given an arbitrarily increasing sequence $(t_n)_{n\in\mathbb N}$ in $(T-\delta,T)$ that converges to $T$, the condition \eqref{t1c}
implies that
\begin{equation}\label{t1d}
e^{-\frac{1}{2}\int_\ell^xb_1(s)\d s}f_1(x)=e^{-\frac{1}{2}\int_\ell^xb_2(s)\d s}f_2(x),\quad x\in(0,\ell)
\end{equation}
and for $b=b_1-b_2$ and $q=q_1-q_2$, we have
\begin{equation}\label{t1ee}
b(x)= b(0)e^{\frac{1}{2}\int_0^{x} (b_1(s)+b_2(s))\d s}+2\int_0^xe^{ \frac{1}{2}\int_\tau^{x} (b_1(s)+b_2(s))\d s} q(\tau)\d\tau,\quad x\in(0,\ell).
\end{equation}
\end{theorem}
\begin{proof} For $j=1,2$, let $A_j=-\partial_x^2 u+V_j$ acting on the space $L^2(0,\ell)$ with domain $D(A_j)=H_0^1(0,\ell)\cap H^2(0,\ell)$. Recall that, the spectrum of $A_j$ consists of an increasing and unbounded sequence $(\lambda_{j,k})_{k\in\mathbb N}$ of simple  eigenvalues and we associate with these eigenvalues a corresponding orthonormal basis $(\phi_{j,k})_{k\in\mathbb N}$ of eigenvalues of $L^2(0,\ell)$. By Proposition \ref{p1}, for $j=1,2$,
the solution of $u_j$ of problem \eqref{eq1}, with  $b=b_j$, $q=q_j$ and $f=f_j$, is lying in $L^2(0,T;H^{2}(0,\ell))$ and, for a.e. $t\in(0,T)$, we have
\begin{equation} \label{form}\partial_xu_j(\ell,t)=\sum_{k=1}^\infty \left(\int_0^t(t-s)^{\alpha-1}E_{\alpha,\alpha}(-\lambda_{j,k}(t-s)^{\alpha})\sigma(s)\d s\right)\langle \tilde{f_j},\phi_{j,k}\rangle\phi_{j,k}'(\ell),
\end{equation}
with $\tilde{f}_j(x) =e^{-\frac{1}{2}\int_\ell^xb_j(s)\d s}f_j(x)$ for $ x\in(0,\ell).$
Using \eqref{t1a}, we get
$$\partial_xu_j(\ell,t)=\sum_{k=1}^\infty \left(\int_0^{T-\delta}(t-s)^{\alpha-1}E_{\alpha,\alpha}(-\lambda_{j,k}(t-s)^{\alpha})\sigma(s)\d s\right)\langle \tilde{f_j},\phi_{j,k}\rangle\phi_{j,k}'(\ell),\quad  t\in(T-\delta,T).$$
From  Proposition \ref{p2}, we deduce that  the extension $h_j$ of $-\partial_xu_j(\ell,t)$ to $\R_+$ given by
$$h_j(t)=\sum_{k=1}^\infty \left(\int_0^{T-\delta}(t-s)^{\alpha-1}E_{\alpha,\alpha}(-\lambda_{j,k}(t-s)^{\alpha})\sigma(s)\d s\right)\langle \tilde{f_j},\phi_{j,k}\rangle\phi_{j,k}'(\ell),\quad t\in(T,+\infty)$$
 is analytic on $(T-\delta,+\infty)$.  Then, condition \eqref{t1b} implies
$h_1(t_n)=h_2(t_n)$ for $n\in\mathbb N.$
Since the sequence $(t_n)_{n\in\mathbb N}$ admits an accumulation point at $T$ and both $h_1$ and $h_2$ are analytic on $(T-\delta,+\infty)$, by fixing $h=h_2-h_1$, we obtain
$$h(t)=0,\quad t\in(T-\delta,+\infty).$$
This proves that the function $h$ is supported on $[0,T-\delta]$ and its Laplace transform $\hat{h}(p)$ is holomorphic with respect to $p\in\mathbb C$. Meanwhile, following Proposition \ref{p2}, one can check
\begin{align*}
\hat{h}(p)=\hat{h}_1(p)-\hat{h}_2(p),\quad p\in\mathbb C_{n_0},
\end{align*}
with $\mathbb C_{n_0}=\{z\in\mathbb C: |z|\not\in \{|\lambda_{j,n}|^{\frac{1}{\alpha}}:\ n\in\mathbb N,\ j=1,2\}\}\setminus (-\infty,0]$ and
\begin{align*}
\hat{h}_j(p)=-\hat{\sigma}(p)\sum_{k=1}^\infty\frac{\langle \tilde{f_j},\phi_{j,k}\rangle\phi_{j,k}'(\ell)}{p^{\alpha}+\lambda_{j,k}},\quad p\in\mathbb C_{n_0}.
\end{align*}
Fix $0<r_1<r_2$, such that $(r_1,r_2)\subset (0,+\infty)\setminus \{|\lambda_{j,n}|^{\frac{1}{\alpha}}:\ n\in\mathbb N,\ j=1,2\} $.
By choosing $p=Re^{\pm i\theta}$, with $R\in (r_1,r_2)$, $\theta\in(0,\pi)$, and sending $\theta\to\pi$, we obtain
\begin{align*}
0&=\hat{h}(-R)-\hat{h}(-R)=\lim_{\theta\to\pi}\hat{h}(Re^{i\theta})-\lim_{\theta\to\pi}\hat{h}(Re^{-i\theta})\\
&=\lim_{\theta\to\pi}\hat{h}_2(Re^{i\theta})-\lim_{\theta\to\pi}\hat{h}_2(Re^{-i\theta})-\left(\lim_{\theta\to\pi}\hat{h}_1(Re^{i\theta})-\lim_{\theta\to\pi}\hat{h}_1(Re^{-i\theta})\right).
\end{align*}
Moreover, we have
\begin{align*}
\lim_{\theta\to\pi}\hat{h}_j(Re^{i\theta})-\lim_{\theta\to\pi}\hat{h}_j(Re^{-i\theta})=2iR^{\alpha}\sin(\alpha\pi)\hat{\sigma}(-R)\sum_{k=1}^\infty\frac{\langle \tilde{f_j},\phi_{j,k}\rangle\phi_{j,k}'(\ell)}{(R^{\alpha}e^{i\alpha\pi}+\lambda_{j,k})(R^{\alpha}e^{-i\alpha\pi}+\lambda_{j,k})}.
\end{align*}
Thus, it follows that
\begin{align}
&R^{\alpha}\sin(\alpha\pi)\hat{\sigma}(-R)\sum_{k=1}^\infty\frac{\langle \tilde{f_1},\phi_{1,k}\rangle\phi_{1,k}'(\ell)}{(R^{\alpha}e^{i\alpha\pi}+\lambda_{1,k})(R^{\alpha}e^{-i\alpha\pi}+\lambda_{1,k})}\nonumber\\
=&R^{\alpha}\sin(\alpha\pi)\hat{\sigma}(-R)\sum_{k=1}^\infty\frac{\langle \tilde{f_2},\phi_{2,k}\rangle\phi_{2,k}'(\ell)}{(R^{\alpha}e^{i\alpha\pi}+\lambda_{2,n})(R^{\alpha}e^{-i\alpha\pi}+\lambda_{2,k})}.
\label{t1e}
\end{align}
Since $\hat{\sigma}$ is holomorphic and non-uniformly vanishing, we may assume without loss of generality that $\hat{\sigma}(-R)\neq0$ for $R\in(r_1,r_2)$. Then, we obtain
\bel{t1g}\sum_{k=1}^\infty\frac{\langle \tilde{f_1},\phi_{1,k}\rangle\phi_{1,k}'(\ell)}{(R^{\alpha}e^{i\alpha\pi}+\lambda_{1,k})(R^{\alpha}e^{-i\alpha\pi}+\lambda_{1,k})}=\sum_{k=1}^\infty\frac{\langle \tilde{f_2},\phi_{2,k}\rangle\phi_{2,k}'(\ell)}{(R^{\alpha}e^{i\alpha\pi}+\lambda_{2,k})(R^{\alpha}e^{-i\alpha\pi}+\lambda_{2,k})},\quad R\in(r_1,r_2).\ee
Meanwhile, for $j=1,2$, since $\tilde{f_j}\in L^2(0,\ell)$, one can check that the map
$$\mathcal G_j:z\mapsto\sum_{k=1}^\infty\frac{\langle \tilde{f_j},\phi_{j,k}\rangle\phi_{j,k}'(\ell)}{(ze^{i\alpha\pi}+\lambda_{j,k})(ze^{-i\alpha\pi}+\lambda_{j,k})}=\partial_x(A_j+ze^{i\alpha\pi})^{-1}(A_j+ze^{-i\alpha\pi})^{-1}\tilde{f_j}(\ell)$$
is meromorphic in $\mathbb C$ with simple poles at $-\lambda_{j,k}e^{\pm i\alpha\pi}$, $k\in\mathbb N$. Therefore, by fixing $\mathcal D=\mathbb C\setminus\{-\lambda_{j,k}e^{\pm i\alpha\pi}:\ j=1,2,\ k\in\mathbb N\}$, we deduce from \eqref{t1g} that
\bel{t1h}\mathcal G_1(z)=\mathcal G_2(z),\quad z\in\mathcal D.\ee
Now consider the set
$$\mathbb K=\{k\in\mathbb N:\ \langle \tilde{f_1},\phi_{1,k}\rangle\neq0\}$$
and observe that since $f_1\not\equiv0$, $\mathbb K\neq\emptyset$.
In light of the identity \eqref{t1h}, for every $k\in \mathbb K$, $-\lambda_{1,k}e^{\pm i\alpha\pi}$ is  a simple pole of $\mathcal G_1$ and by \eqref{t1h} it is also a simple pole of $\mathcal G_2$. Therefore, we have
$$\{\lambda_{1,k}e^{\pm i\alpha\pi}:\ k\in\mathbb K\}\subset \{\lambda_{2,k}e^{\pm i\alpha\pi}:\ k\in\mathbb N\},$$
which implies
$$\{\lambda_{1,k}:\ k\in\mathbb K\}\subset \{\lambda_{2,k}:\ k\in\mathbb N\}.$$
Therefore, we can find  a map $\psi$, bijective from  $\mathbb N$ to $\mathbb N$, such that
$$\lambda_{1,k}=\lambda_{2,\psi(k)},\quad k\in\mathbb K.$$
For each $\mu>0$, let
$m(\mu):=\sharp \{k\in\mathbb N:\ \lambda_{1,k}\leq \mu\ \textrm{and }\langle \tilde{f_1},\phi_{1,k}\rangle=0\}$.
Now we define
$$N_j(\mu)=\sharp \{k\in\mathbb N:\ \lambda_{j,k}\leq \mu\},\quad j=1,2,\quad \mbox{and}\quad S(\mu)=\sharp \{k\in\mathbb N:\ \lambda_{1,k}\leq \mu,\ \lambda_{1,k}=\lambda_{2,\psi(k)}\}.$$
In view of \eqref{l1a} in Lemma \ref{l1} in the appendix, we have
$$N_j(\mu)=\frac{\ell\sqrt{\mu}}{\pi}+\underset{\mu\to+\infty}{ \mathcal O}(1),\quad j=1,2.$$
By combining this asymptotic with \eqref{l1b} from Lemma \ref{l1} in the appendix, we deduce that there exists $\mu_0>0$ such that for all $\mu>\mu_0$ and any $\varepsilon \in (0,\frac{\ell}{2r}-\frac{1}{2})$, we have
\begin{align*}
	S(\mu)&= N_1(\mu)-m(\mu)\geq N_1(\mu) - (1+2\varepsilon)\frac{r}{\pi}\sqrt{\mu}\\
	&\geq (1-(1+\varepsilon)\frac{r}{\ell})N_j(\mu)+\frac{(1+\varepsilon)}{2}\frac{r}{\ell},\quad j=1,2.
\end{align*}
Thus, in view of condition \eqref{t1b},  by applying Proposition \ref{p4} and choosing $\varepsilon$  sufficiently small, we deduce that $V_1=V_2$, i.e., $A_1=A_2=A$. Thus we can consider $(\phi_k)_{ k\in\mathbb N}$ that forms an orthonormal basis in $L^2(0,\ell)$ of eigenfunctions  associated with the eigenvalues $(\lambda_k)_{k\in\mathbb N}$ of $A$ such that \eqref{t1h} holds true for $\phi_{1,k}=\phi_{2,k}=\phi_k$ and $\lambda_k=\lambda_{1,k}=\lambda_{2,k}$. Then, the identity \eqref{t1h} can be rewritten as
$$\sum_{k=1}^\infty\frac{\langle \tilde{f_1},\phi_{k}\rangle\phi_{k}'(\ell)}{(ze^{i\alpha\pi}+\lambda_{k})(ze^{-i\alpha\pi}+\lambda_{k})}=\sum_{k=1}^\infty\frac{\langle \tilde{f_2},\phi_{k}\rangle\phi_{k}'(\ell)}{(ze^{i\alpha\pi}+\lambda_{k})(ze^{-i\alpha\pi}+\lambda_{k})},\quad z\in\mathbb C\setminus\{-\lambda_ke^{\pm i\alpha\pi}:\  k\in\mathbb N\}.$$
Using this last identity,  we can prove by iteration that
$$\langle \tilde{f}_1,\phi_k\rangle=\langle \tilde{f}_2,\phi_k\rangle,\quad k\in\mathbb N$$
and it follows that
$$e^{-\frac{1}{2}\int_\ell^xb_1(s)\d s}f_1(x)=\tilde{f}_1(x)=\tilde{f}_2(x)=e^{-\frac{1}{2}\int_\ell^xb_2(s)\d s}f_2(x),\quad x\in(0,\ell).$$
This proves the identity \eqref{t1d}. Finally, using the fact that
$$-\tfrac{1}{2}b_1'(x)+\tfrac{1}{4}b_1^2(x)+q_1(x)=V_1(x)=V_2(x)=-\tfrac{1}{2}b_2'(x)+\tfrac{1}{4}b_2^2(x)+q_2(x),\quad x\in(0,\ell)$$
we deduce that $b=b_1-b_2$ solves the ODE
$$b'(x)-\tfrac{1}{2}(b_1(x)+b_2(x))b=2(q_1(x)-q_2(x)),\quad x\in(0,\ell)$$
whose solutions take the form \eqref{t1ee}. This completes the proof of the theorem.
\end{proof}
\subsection{Proof of Theorem \ref{t2_intro}}
We prove a slightly stronger version of Theorem~\ref{t2_intro} as follows.
\begin{theorem}\label{t2} Let the condition of Theorem \ref{t1} hold with $\alpha\in(0,1)\setminus \mathbb Q$.
Consider $\sigma_j\in L^2(0,T)$, $j=1,2$, be a non-uniformly vanishing function such that  \eqref{t1a} is fulfilled with $\sigma=\sigma_j$ and denote by $u_j$ the solution of \eqref{eq1}, with  $\sigma=\sigma_j$, $b=b_j$, $f=f_j$ and $q=q_j$. Then condition \eqref{t1c} implies that there exists $\beta\in\R\setminus\{0\}$ such that
there hold \eqref{t1ee}  and
\begin{align}
\label{t2a}\sigma_2(t)&=\beta\sigma_1(t),\quad t\in(0,T),\\
\label{t2b} e^{-\frac{1}{2}\int_\ell^xb_1(s)\d s}f_1(x)&=\beta e^{-\frac{1}{2}\int_\ell^xb_2(s)\d s}f_2(x),\quad x\in(0,\ell).
\end{align}
\end{theorem}
\begin{proof}
We use the same notations as in Section \ref{ssec:thm1}. By repeating the argument of Theorem \ref{t1}, for all $R \in(0,+\infty)\setminus \{|\lambda_{j,n}|^{\frac{1}{\alpha}}:\ n\in\mathbb N,\ j=1,2\}$, we get
\begin{align*}
&R^{\alpha}\sin(\alpha\pi)\hat{\sigma_1}(-R)\sum_{k=1}^\infty\frac{\langle \tilde{f_1},\phi_{1,k}\rangle\phi_{1,k}'(\ell)}{(R^{\alpha}e^{i\alpha\pi}+\lambda_{1,k})(R^{\alpha}e^{-i\alpha\pi}+\lambda_{1,k})}\\
=&R^{\alpha}\sin(\alpha\pi)\hat{\sigma_2}(-R)\sum_{k=1}^\infty\frac{\langle \tilde{f_2},\phi_{2,k}\rangle\phi_{2,k}'(\ell)}{(R^{\alpha}e^{i\alpha\pi}+\lambda_{2,n})(R^{\alpha}e^{-i\alpha\pi}+\lambda_{2,k})}.
\end{align*}
 For all $r\in(0,+\infty)\setminus \{|\lambda_{j,n}|:\ n\in\mathbb N,\ j=1,2\}$, by fixing $R=r^{1\over\alpha}$, we find
\bel{t2d}\hat{\sigma}_1(-r^{1\over\alpha})\sum_{k=1}^\infty\frac{\langle \tilde{f_1},\phi_{1,k}\rangle\phi_{1,k}'(\ell)}{(re^{i\alpha\pi}+\lambda_{1,k})(re^{-i\alpha\pi}+\lambda_{1,k})}-
\hat{\sigma}_2(-r^{1\over\alpha})\sum_{k=1}^\infty\frac{\langle \tilde{f_2},\phi_{2,k}\rangle\phi_{2,k}'(\ell)}{(re^{i\alpha\pi}+\lambda_{2,n})(re^{-i\alpha\pi}+\lambda_{2,k})}\equiv 0.\ee
Let $\mathcal U:=\mathbb C\setminus( \{-e^{ \pm i\alpha\pi}\lambda_{j,k}:\ k\in\mathbb N,\ j=1,2\}\cup (-\infty,0])$ and consider
the
map
$$\mathcal G:\mathcal U\ni
 z\mapsto \hat{\sigma}_1(-z^{1\over\alpha})\Psi_1(\cdot,z)+
\hat{\sigma}_2(-z^{1\over\alpha})\Psi_2(\cdot,z),$$
with
$$
\Psi_1(\cdot,z)=\sum_{k=1}^\infty\frac{\langle \tilde{f_1},\phi_{1,k}\rangle\phi_{1,k}'(\ell)}{(ze^{i\alpha\pi}+\lambda_{1,k})(ze^{-i\alpha\pi}+\lambda_{1,k})}\quad \mbox{and}\quad
\Psi_2(\cdot,z)=-\sum_{k=1}^\infty\frac{\langle \tilde{f_2},\phi_{2,k}\rangle\phi_{2,k}'(\ell)}{(ze^{i\alpha\pi}+\lambda_{2,n})(ze^{-i\alpha\pi}+\lambda_{2,k})}.
$$
Repeating the argument of Theorem \ref{t1} gives that the map $\mathcal G$ is  holomorphic in $\mathcal U$. Thus, condition \eqref{t2d} and unique continuation of holomorphic functions imply that $\mathcal G$ vanishes identically in $\mathcal U$. By choosing $z=\varrho e^{\pm i\theta}$   and sending $\theta$ to $\pi$, we get
\begin{align}\label{t2e}
\hat{\sigma}_1(-\varrho^{1\over\alpha}e^{i\pi\over\alpha})\Psi_1(\cdot,-\varrho)+
\hat{\sigma}_2(-\varrho^{1\over\alpha}e^{i\pi\over\alpha})\Psi_2(\cdot,-\varrho)&\equiv 0,
\\ \label{t2f}
\hat{\sigma}_1(-\varrho^{1\over\alpha}e^{-{i\pi\over\alpha}})\Psi_1(\cdot,-\varrho)+
\hat{\sigma}_2(-\varrho^{1\over\alpha}e^{-i{\pi\over\alpha}})\Psi_2(\cdot,-\varrho)&\equiv 0,
\end{align}
where $\varrho\in \{s>0:\ s\neq \lambda_{j,k},\ k\in\mathbb N,\ j=1,2\}$. By the continuity, the identities \eqref{t2e}-\eqref{t2f} can then be extended to any $\varrho>0$. Note that $\Psi_j\not\equiv0$, $j=1,2$, since otherwise applying the argument in the last step of Theorem \ref{t1} this would imply that $f_1=f_2\equiv0$, contradicting the assumptions. By combining this with the fact that, for $j=1,2$,  the map $\Psi_j$ is holomorphic in
  the set $\mathcal U_1=\C\setminus  \{-e^{ \pm i\alpha\pi}\lambda_{j,k}:\ k\in\mathbb N,\ j=1,2\}$, we deduce that
	there exists   $\rho_j>0$ such that $|\Phi_{j}(\cdot,z)|> 0$ for $0<|z|<\rho_j$. Then, the determinant of the system \eqref{t2e}, \eqref{t2f}
  is zero for $0<\varrho<r_0=\min(\rho_1,\rho_2)$, i.e.,
\bea
\label{t2g}  \left|\left[\begin{array}{lll}
&\hat{\sigma}_1(-\varrho^{1\over\alpha}e^{i\pi\over\alpha}) &\hat{\sigma}_2(-\varrho^{1\over\alpha}e^{i\pi\over\alpha})
\\
&\hat{\sigma}_1(-\varrho^{1\over\alpha}e^{-{i\pi\over\alpha}}) &\hat{\sigma}_2(-\varrho^{1\over\alpha}e^{-{i\pi\over\alpha}})
\end{array}\right]\right|= 0,\quad 0<\varrho<r_0.
\eea
Since $\sigma_j$, $j=1,2$, are compactly supported and non-uniformly vanishing, $\hat{\sigma}_j$ are necessarily entire and not uniformly vanishing. Moreover, we can find $r_1\in(0,r_0)$ such that
$$\min(|\hat{\sigma}_1(z)|,|\hat{\sigma}_2(z)|)>0,\quad 0<|z|<r_1.$$
Thus, either $\hat{\sigma}_1(z)/\hat{\sigma}_2(z)$ or $\hat{\sigma}_2(z)/\hat{\sigma}_1(z)$ will be bounded as $z\to0$ and then can be extended to a holomorphic function with respect to  $z\in D_{r_1}:=\{\eta\in\mathbb C:\ |\eta|<r_1\}$. Without loss of generality, we may assume that $\hat{\sigma}_2(z)/\hat{\sigma}_1(z)$ admits a holomorphic extension with respect to  $z\in D_{r_1}$ into a map denoted by $h$.
Then, the identity \eqref{t2g} can be rewritten as
$$
\hat{\sigma}_1(-\varrho^{1\over\alpha}e^{i\pi\over\alpha})\hat{\sigma}_1(-\varrho^{1\over\alpha}e^{-{i\pi\over\alpha}})
\left(h(-\varrho^{1\over\alpha}e^{i\pi\over\alpha})-h(-\varrho^{1\over\alpha}e^{-{i\pi\over\alpha}})\right)= 0,
\quad 0<\varrho<r_1.
$$
This implies that $h(-\varrho^{1\over\alpha}e^{i\pi\over\alpha})= h(-\varrho^{1\over\alpha}e^{-{i\pi\over\alpha}})$ for $0<\varrho<r_1$. Then, since the map $h$ is holomorphic in  $D_{r_1}$, we get
$
h(z)=h(ze^{-{2i\pi\over\alpha}})$ for $    z\in D_{r_1^{1\over\alpha}}.
$
By differentiating the relation at $z=0$, we have
$
h^{(k)}(0)=h^{(k)}(0)e^{-{2ki\pi\over\alpha}}$ for $k\in\mathbb N.$
Since $\alpha$ is irrational, we have $e^{-{2ki\pi\over\alpha}}\ne 1$, $k\in\mathbb N$, which implies  $h^{(k)}(0)=0$, $k\in\mathbb N$. Hence, the function $h$ is a constant map which will necessary be real values and non-vanishing which clearly implies \eqref{t2a}. By replacing $f_2$ by $\beta f_2$ and applying Theorem \ref{t1} we obtain \eqref{t2b} and \eqref{t1ee}. This completes the proof of the theorem.\end{proof}

\subsection{Proof of Theorem \ref{tt1_intro}}

We give the proof of a stronger version than  Theorem~\ref{tt1_intro} stated  as follows.

\begin{theorem}\label{tt1} For $j=1,2$, let  $b_j\in C^1([0,\ell])$, $q_j\in C([0,\ell])$. Consider $f_j\in L^2(0,\ell)$, $\sigma_j\in L^2(0,T)$, $j=1,2$, be a non-uniformly vanishing function such that  \eqref{t1a} is fulfilled with $\sigma=\sigma_j$ and denote by $u_j$ the solution of \eqref{eq1}, with  $\sigma=\sigma_j$, $b=b_j$, $f=f_j$ and $q=q_j$. Moreover, if $\alpha\in(0,1)\cap \mathbb Q$ then we assume that $\sigma_1=\sigma_2$.
We fix also $V_j=-\frac{1}{2}b_j'+\frac{1}{4}b_j^2+q_j$, $a_1,a_2\in(0,\ell)$, with $a_1<a_2$, $a_1< \frac{\ell}{2}$ and $a_1+a_2\geq\ell$,  and we assume that
\bel{tt1aa}V_1(x)=V_2(x),\quad x\in[a_1,a_2], \quad b_1(a_2)=b_2(a_2).\ee
Also assume that there exists $r\in (0,\min(\ell-2a_1,2a_2-\ell))$ such that $\mathrm{supp}(f_1)\subset [0,r]$.
Then,  given any increasing sequence $(t_n)_{n\in\mathbb N}$ in $(T-\delta,T)$ that converges to $T$ and any $\epsilon\in (0,\min(a_2,\ell-a_2))$, the condition \eqref{tt1a}
implies  that there exists $\beta\in\R\setminus\{0\}$ such that the conditions \eqref{t1ee},
\eqref{t2a} and
\bel{tt1ab} e^{-\frac{1}{2}\int_{a_2}^xb_1(s)\d s}f_1(x)=\beta e^{-\frac{1}{2}\int_{a_2}^xb_2(s)\d s}f_2(x),\quad x\in(0,\ell)\ee
 are fulfilled.
\end{theorem}

\begin{proof}

Consider first the case $\alpha\in\mathbb Q\cap (0,1)$ and $\sigma_1=\sigma_2=\sigma$. We use the notation in the proof of Theorem \ref{t1}. In view of Proposition \ref{p1}, for $j=1,2$,
the solution $u_j$ of problem \eqref{eq1} with  $b=b_j$, $q=q_j$ and $f=f_j$ belongs to $L^2(0,T;H^{2}(0,\ell))$ and, for a.e. $t\in(0,T)$, we have
\bel{tt1d}u_j(x,t)=e^{\frac{1}{2}\int_\ell^xb_j(s)\d s}\sum_{k=1}^\infty \left(\int_0^t(t-s)^{\alpha-1}E_{\alpha,\alpha}(-\lambda_{j,k}(t-s)^{\alpha})\sigma(s)\d s\right)\langle \tilde{f_j},\phi_{j,k}\rangle\phi_{j,k}(x),\quad   x\in(0,\ell).\ee
Using \eqref{t1a}, for all $t\in(T-\delta,T)$ and all $ x\in[a_2-\epsilon,a_2+\epsilon]$, we obtain
\begin{align*}
u_j(x,t)=e^{\frac{1}{2}\int_\ell^xb_j(s)\d s}\sum_{k=1}^\infty \left(\int_0^{T-\delta}(t-s)^{\alpha-1}E_{\alpha,\alpha}(-\lambda_{j,k}(t-s)^{\alpha})\sigma(s)\d s\right)\langle \tilde{f_j},\phi_{j,k}\rangle\phi_{j,k}(x).
\end{align*}
From Proposition \ref{p2}, we deduce that  the extension $w_j$ of $u_j|_{[a_2-\epsilon,a_2+\epsilon]\times(0,T)}$ to $[a_2-\epsilon,a_2+\epsilon]\times\R_+$ given, for all $t\in(T,+\infty)$ and all $ x\in[a_2-\epsilon,a_2+\epsilon]$, by
$$w_j(x,t)=e^{\frac{1}{2}\int_\ell^xb_j(s)\d s}\sum_{k=1}^\infty \left(\int_0^{T-\delta}(t-s)^{\alpha-1}E_{\alpha,\alpha}(-\lambda_{j,k}(t-s)^{\alpha})\sigma(s)\d s\right)\langle \tilde{f_j},\phi_{j,k}\rangle\phi_{j,k}(x),$$
 is analytic with respect to $t\in(T-\delta,+\infty)$ as a  map taking values in $C([a_2-\epsilon,a_2+\epsilon])$.  Then, condition \eqref{tt1a} implies
$$w_1(x,t_n)=w_2(x,t_n),\quad n\in\mathbb N,\ x\in[a_2-\epsilon,a_2+\epsilon].$$
Since the sequence $(t_n)_{n\in\mathbb N}$ admits an accumulation point at $T$ and both $w_1$ and $w_2$ are analytic  with respect to $t\in(T-\delta,+\infty)$, by fixing $w=w_2-w_1$, we obtain
$$w(x,t)=0,\quad t\in(T-\delta,+\infty),\ x\in[a_2-\epsilon,a_2+\epsilon].$$
This proves that the function $w$ is supported on $[a_2-\epsilon,a_2+\epsilon]\times[0,T-\delta]$ and, for all $x\in[a_2-\epsilon,a_2+\epsilon]$, its Laplace transform in time $\hat{w}(x,p)=\int_0^{+\infty}w(x,t)e^{-pt}\d t$ is holomorphic with respect to $p\in\mathbb C$. Meanwhile, following Proposition \ref{p2}, one can check
$$\hat{w}(x,p)=\hat{w}_1(x,p)-\hat{w}_2(x,p),\quad p\in\mathbb C_{n_0},\ x\in[a_2-\epsilon,a_2+\epsilon],$$
with
$$\hat{w}_j(x,p)=-\hat{\sigma}(p)e^{\frac{1}{2}\int_\ell^xb_j(s)\d s}\sum_{k=1}^\infty\frac{\langle \tilde{f_j},\phi_{j,k}\rangle\phi_{j,k}(x)}{p^{\alpha}+\lambda_{j,k}},\quad p\in\mathbb C_{n_0},\ x\in(0,\ell),$$
where $\mathbb C_{n_0}=\{z\in\mathbb C: |z|\not\in \{|\lambda_{j,n}|^{\frac{1}{\alpha}}:\ n\in\mathbb N,\ j=1,2\}\}\setminus (-\infty,0]$.
By choosing $p=Re^{\pm i\theta}$, with $R \in(0,+\infty)\setminus \{|\lambda_{j,n}|^{\frac{1}{\alpha}}:\ n\in\mathbb N,\ j=1,2\}$, $\theta\in(0,\pi)$, and sending $\theta\to\pi$, we obtain
$$\begin{aligned}0=&\hat{w}(x,-R)-\hat{w}(x,-R)
=\lim_{\theta\to\pi}\hat{w}_2(x,Re^{i\theta})-\lim_{\theta\to\pi}\hat{w}_2(x,Re^{-i\theta})\\
&-(\lim_{\theta\to\pi}\hat{w}_1(x,Re^{i\theta})-\lim_{\theta\to\pi}\hat{w}_1(x,Re^{-i\theta})),\ x\in[a_2-\epsilon,a_2+\epsilon].\end{aligned}$$
Moreover, we have
\begin{align*}
&\lim_{\theta\to\pi}\hat{w}_j(x,Re^{i\theta})-\lim_{\theta\to\pi}\hat{w}_j(x,Re^{-i\theta})\\
=&2iR^{\alpha}\sin(\alpha\pi)\hat{\sigma}(-R)e^{\frac{1}{2}\int_\ell^xb_j(s)\d s}\sum_{k=1}^\infty\frac{\langle \tilde{f_j},\phi_{j,k}\rangle\phi_{j,k}(x)}{(R^{\alpha}e^{i\alpha\pi}+\lambda_{j,k})(R^{\alpha}e^{-i\alpha\pi}+\lambda_{j,k})},\ x\in(0,\ell),
\end{align*}
and it follows that for $x\in[a_2-\epsilon,a_2+\epsilon]$,
\begin{align}&R^{\alpha}\sin(\alpha\pi)\hat{\sigma}(-R)e^{\frac{1}{2}\int_\ell^xb_1(s)\d s}\sum_{k=1}^\infty\frac{\langle \tilde{f_1},\phi_{1,k}\rangle\phi_{1,k}(x)}{(R^{\alpha}e^{i\alpha\pi}+\lambda_{1,k})(R^{\alpha}e^{-i\alpha\pi}+\lambda_{1,k})}\nonumber\\
=&R^{\alpha}\sin(\alpha\pi)\hat{\sigma}(-R)e^{\frac{1}{2}\int_\ell^xb_2(s)\d s}\sum_{k=1}^\infty\frac{\langle \tilde{f_2},\phi_{2,k}\rangle\phi_{2,k}(x)}{(R^{\alpha}e^{i\alpha\pi}+\lambda_{2,n})(R^{\alpha}e^{-i\alpha\pi}+\lambda_{2,k})}.\label{tt1e}
\end{align}
Since $\hat{\sigma}$ is holomorphic and non-uniformely vanishing, we can find $r_1<r_2$ such that $\hat{\sigma}(-R)\neq0$ for $R\in(r_1,r_2)$ and $(r_1,r_2)\subset (0,+\infty)\setminus \{|\lambda_{j,n}|^{\frac{1}{\alpha}}:\ n\in\mathbb N,\ j=1,2\}$. Thus, for all $x\in[a_2-\epsilon,a_2+\epsilon]$, we obtain
$$\begin{aligned}&e^{\frac{1}{2}\int_\ell^xb_1(s)\d s}\sum_{k=1}^\infty\frac{\langle \tilde{f_1},\phi_{1,k}\rangle\phi_{1,k}(x)}{(R^{\alpha}e^{i\alpha\pi}+\lambda_{1,k})(R^{\alpha}e^{-i\alpha\pi}+\lambda_{1,k})}\\
&=e^{\frac{1}{2}\int_\ell^xb_2(s)\d s}\sum_{k=1}^\infty\frac{\langle \tilde{f_2},\phi_{2,k}\rangle\phi_{2,k}(x)}{(R^{\alpha}e^{i\alpha\pi}+\lambda_{2,k})(R^{\alpha}e^{-i\alpha\pi}+\lambda_{2,k})},\quad R\in(r_1,r_2).\end{aligned}$$
Further, for $j=1,2$,  since $\tilde{f_j}\in L^2(0,\ell)$, one can check that the map
$$\mathcal H_j:z\mapsto\sum_{k=1}^\infty\frac{\langle \tilde{f_j},\phi_{j,k}\rangle\phi_{j,k}|_{[a_2-\epsilon,a_2+\epsilon]}}{(ze^{i\alpha\pi}+\lambda_{j,k})(ze^{-i\alpha\pi}+\lambda_{j,k})}$$
is meromorphic in $\mathbb C$, as a map taking values in $C([a_2-\epsilon,a_2+\epsilon])$, with simple poles at $-\lambda_{j,k}e^{\pm i\alpha\pi}$, $k\in\mathbb N$. Therefore, in a similar way to Theorem \ref{t1}, we get
\bel{tt1h}e^{\frac{1}{2}\int_\ell^xb_1(s)\d s}\mathcal H_1(z,x)=e^{\frac{1}{2}\int_\ell^xb_2(s)\d s}\mathcal H_2(z,x),\quad z\in\mathcal D,\ x\in[a_2-\epsilon,a_2+\epsilon].\ee
By following the argument in Theorem \ref{t1} and fixing $\mathbb K=\{k\in\mathbb N:\ \langle \tilde{f_1},\phi_{1,k}\rangle\neq0\}$,
we observe that there exists a map $\psi$, bijective from  $\mathbb N$ to $\mathbb N$, such that, for all $x\in[a_2-\epsilon,a_2+\epsilon]$, we have
$$\lambda_{1,n}=\lambda_{2,\psi(n)},\  e^{\frac{1}{2}\int_\ell^xb_1(s)\d s}\langle \tilde{f_1},\phi_{1,n}\rangle\phi_{1,n}(x)=e^{\frac{1}{2}\int_\ell^xb_2(s)\d s}\langle \tilde{f_2},\phi_{2,\psi(n)}\rangle\phi_{2,\psi(n)}(x),\quad  n\in\mathbb K.$$
Therefore, by applying \eqref{tt1aa}, we get
\begin{align*}
e^{\frac{1}{2}\int_\ell^{a_2}b_1(s)\d s}\langle \tilde{f_1},\phi_{1,n}\rangle\phi_{1,n}(a_2)&=e^{\frac{1}{2}\int_\ell^{a_2}b_2(s)\d s}\langle \tilde{f_2},\phi_{2,\psi(n)}\rangle\phi_{2,\psi(n)}(a_2),\\
e^{\frac{1}{2}\int_\ell^{a_2}b_1(s)\d s}\langle \tilde{f_1},\phi_{1,n}\rangle\phi_{1,n}'(a_2)&=e^{\frac{1}{2}\int_\ell^{a_2}b_2(s)\d s}\langle \tilde{f_2},\phi_{2,\psi(n)}\rangle\phi_{2,\psi(n)}'(a_2),
\end{align*}
and it follows that
\bel{tt1f}\lambda_{1,n}=\lambda_{2,\psi(n)},\quad  \phi_{1,n}(a_2)\phi_{2,\psi(n)}'(a_2)=\phi_{2,\psi(n)}(a_2)\phi_{1,n}'(a_2),\quad n\in\mathbb K.\ee
Then, in a similar way to Theorem \ref{t1}, by fixing
$$E(a_2,\mu)=\sharp \{k\in\mathbb N:\ \lambda_{1,k}\leq \mu,\ \lambda_{1,k}=\lambda_{2,\psi(k)},\ \phi_{1,k}(a_2)\phi_{2,\psi(k)}'(a_2)=\phi_{2,\psi(k)}(a_2)\phi_{1,k}'(a_2)\},$$
$$N_j(\mu)=\sharp \{k\in\mathbb N:\ \lambda_{j,k}\leq \mu\},\quad j=1,2,$$
$$S(\mu)=\sharp \{k\in\mathbb N:\ \lambda_{1,k}\leq \mu,\ \lambda_{1,k}=\lambda_{2,\psi(k)}\},\ m(\mu):=\sharp \{k\in\mathbb N:\ \lambda_{1,k}\leq \mu\ \textrm{and }\langle \tilde{f_1},\phi_{1,k}\rangle=0\},$$
we can find $\mu_0>0$ such that for all $\mu>\mu_0$ and any $\varepsilon \in (0,\frac{\ell}{2r}-\frac{1}{2})$, we have
\begin{align*}
	\min(S(\mu),E(a_2,\mu))&= N_1(\mu)-m(\mu)\geq N_1(\mu) - (1+2\varepsilon)\frac{r}{\pi}\sqrt{\mu}\\
	&\geq (1-(1+\varepsilon)\frac{r}{\ell})N_j(\mu)+\frac{(1+\varepsilon)}{2}\frac{r}{\ell},\quad j=1,2.
\end{align*}
In light of condition \eqref{tt1aa} and noting $r<\min(\ell-2a_1,2a_2-\ell)$, by choosing $\varepsilon$ to be sufficiently small, we get
$$(1-(1+\varepsilon)\frac{r}{\ell})\geq \max(2\frac{a_1}{\ell},2(1-\frac{a_2}{\ell})).$$
Thus, conditions \eqref{p5a} and \eqref{p5b} of the appendix are fulfilled and Proposition \ref{p5} implies  $V_1=V_2$. Then, repeating the argument of Theorem \ref{t1}, for all $x\in[a_2-\epsilon,a_2+\epsilon]$ and all $z\in\mathbb C\setminus\{-\lambda_ke^{\pm i\alpha\pi}:\  k\in\mathbb N\}$,  the identity \eqref{tt1h} can be rewritten as
$$e^{\frac{1}{2}\int_\ell^xb_1(s)\d s}\sum_{k=1}^\infty\frac{\langle \tilde{f_1},\phi_{1,k}\rangle\phi_{k}(x)}{(ze^{i\alpha\pi}+\lambda_{1,k})(ze^{-i\alpha\pi}+\lambda_{1,k})}=e^{\frac{1}{2}\int_\ell^xb_2(s)\d s}\sum_{k=1}^\infty\frac{\langle \tilde{f_2},\phi_{k}\rangle\phi_{1,k}(x)}{(ze^{i\alpha\pi}+\lambda_{1,k})(ze^{-i\alpha\pi}+\lambda_{1,k})}.$$
From this last identity,  we get
$$e^{\frac{1}{2}\int_\ell^{x}b_1(s)\d s}\langle \tilde{f}_1,\phi_{1,k}\rangle\phi_{1,k}(x)=e^{\frac{1}{2}\int_\ell^{x}b_2(s)\d s}\langle \tilde{f}_2,\phi_{1,k}\rangle\phi_{1,k}(x),\quad k\in\mathbb N, \ x\in[a_2-\epsilon,a_2+\epsilon],$$
and, since for all $k\in\mathbb N$ and all $\delta\in(0,\epsilon)$, $\phi_{1,k}|_{[a_2-\delta,a_2+\delta]}\not\equiv0$, we obtain
\begin{align*}
e^{\frac{1}{2}\int_\ell^{a_2}b_1(s)\d s}\langle \tilde{f}_1,\phi_{1,k}\rangle=e^{\frac{1}{2}\int_\ell^{a_2}b_2(s)\d s}\langle \tilde{f}_2,\phi_{1,k}\rangle,\quad k\in\mathbb N.
\end{align*}
Consequently, it follows that
$$e^{-\frac{1}{2}\int_{a_2}^xb_1(s)\d s}f_1(x)=e^{\frac{1}{2}\int_\ell^{a_2}b_1(s)\d s}\tilde{f}_1(x)=e^{\frac{1}{2}\int_\ell^{a_2}b_2(s)\d s}\tilde{f}_2(x)=e^{-\frac{1}{2}\int_{a_2}^xb_2(s)\d s}f_2(x),\quad x\in(0,\ell)$$
and we deduce that  \eqref{t2a},  \eqref{t1ee} and \eqref{tt1ab} are fulfilled with $\beta=1$.

Now consider the case $\alpha\in(0,1)\setminus \mathbb Q$. In this case, for all $x\in[a_2-\epsilon,a_2+\epsilon]$ and all $R\in(0,+\infty)\setminus \{|\lambda_{j,n}|^{\frac{1}{\alpha}}:\ n\in\mathbb N,\ j=1,2\}\}$, we obtain
\begin{align*}&R^{\alpha}\sin(\alpha\pi)\hat{\sigma_1}(-R)e^{\frac{1}{2}\int_\ell^{x}b_1(s)\d s}\sum_{k=1}^\infty\frac{\langle \tilde{f_1},\phi_{1,k}\rangle\phi_{1,k}(x)}{(R^{\alpha}e^{i\alpha\pi}+\lambda_{1,k})(R^{\alpha}e^{-i\alpha\pi}+\lambda_{1,k})}\\
=&R^{\alpha}\sin(\alpha\pi)\hat{\sigma_2}(-R)e^{\frac{1}{2}\int_\ell^{x}b_2(s)\d s}\sum_{k=1}^\infty\frac{\langle \tilde{f_2},\phi_{2,k}\rangle\phi_{2,k}(x)}{(R^{\alpha}e^{i\alpha\pi}+\lambda_{2,n})(R^{\alpha}e^{-i\alpha\pi}+\lambda_{2,k})}.\end{align*}
Since for all $k\in\mathbb N$ and $j=1,2$, $\phi_{j,k}|_{[a_2-\epsilon,a_2+\epsilon]}\not\equiv0$ and following the argument for Theorem \ref{t2}, we can first prove that there exists $\beta\in\R\setminus\{0\}$ such that \eqref{t2a} is fulfilled. Then, replacing $f_2$ by $\beta f_2$ and repeating the arguments for the case $\alpha\in\mathbb Q\cap (0,1)$ and $\sigma_1=\sigma_2=\sigma$, we deduce that \eqref{t2a}, \eqref{t1ee} and \eqref{tt1ab} are fulfilled.\end{proof}

\subsection{Proof of Theorem \ref{t4_intro}}
We prove a slightly stronger version as follows.
\begin{theorem}\label{t4} For $j=1,2$,  $b_j\in C^1([0,\ell])$, $q_j\in C([0,\ell])$. Consider $g_j\in L^2(\Omega)$ and $\sigma_j\in L^2(0,T)$, $j=1,2$,  non-uniformly vanishing functions such that  there exists $\delta\in(0,T)$ satisfying \eqref{t1a} with $\sigma=\sigma_j$. In addition, if $\alpha\in(0,1)\cap \mathbb Q$ then we assume that $\sigma_1=\sigma_2$. We fix also $V_j=-\frac{1}{2}b_j'+\frac{1}{4}b_j^2+q_j$ and we assume that there exists $r\in(0,\ell)$ such that supp$(f_1)\subset [0,r]\times\omega$ and \eqref{t1b} is fulfilled. Denote by $v_j$ the solution of \eqref{eq111}, with  $\sigma=\sigma_j$, $b=b_j$, $g=g_j$ and $q=q_j$. Then the condition
\bel{t4a}\partial_{x_1}v_1(\ell,x',t_n)=\partial_{x_1}v_2(\ell,x',t_n),\quad n\in\mathbb N,\ x'\in\omega\ee
implies  that there exists a constant $\beta\in\R\setminus\{0\}$ such that \eqref{t1ee},
\eqref{t2a} and
\bel{t4b} e^{-\frac{1}{2}\int_\ell^{x_1}b_1(s)\d s}g_1(x_1,x')=\beta e^{-\frac{1}{2}\int_\ell^{x_1}b_2(s)\d s}g_2(x_1,x'),\quad x=(x_1,x')\in\Omega\ee
are fulfilled.
\end{theorem}
\begin{proof} Let $\Delta'=\partial_{x_2}^2+\ldots+\partial_{x_n}^2$ be the Laplacian in $\omega$ and let $H$ be the unbounded operator acting on $L^2(\omega)$ with its domain $D(H)=H^1_0(\omega)\cap H^2(\omega)$ defined by $Hg=-\Delta'g,\quad g\in D(H).$ Then $H$ is a selfadjoint operator whose  spectrum
consists of an increasing sequence of strictly positive  eigenvalues
$(\mu_{m})_{m\geq1}$ of finite multiplicity.    In the Hilbert space $L^2(\omega)$, for each eigenvalue $\mu_m$, we fix also $\ell_m\in\mathbb N$   the algebraic multiplicity of $\mu_m$ and the family $\{\psi_{m,k}\}_{k=1}^{\ell_m}$ of eigenfunctions of $H$,
which forms an orthonormal basis in $L^2(\omega)$ of the algebraic eigenspace of $H$ associated with $\mu_m$. Then $\{\psi_{m,k}:\ m\in\mathbb N,\ k=1,\ldots,\ell_m\}$ forms an orthonormal basis in $L^2(\omega)$. For all $m\in\mathbb N$, $k=1,\ldots,\ell_m$ and $j=1,2$, we denote by $f_{k,m}^j$ the map defined on $(0,\ell)$ by
\begin{equation*}
f_{m,k}^j:(0,\ell)\ni s\mapsto \left\langle g_j(s,\cdot),\psi_{m,k}\right\rangle_{L^2(\omega)}. 
\end{equation*}
In the same way, we define the map $u_{m,k}^j$  on $(0,\ell)\times(0,T)$ by
\begin{equation*}
u_{m,k}^j:(0,\ell)\times(0,T)\ni (s,t)\mapsto \left\langle v_j(s,\cdot,t),\psi_{m,k}\right\rangle_{L^2(\omega)}. 
\end{equation*}
One can easily check that, for all $m\in\mathbb N$ and $j=1,2$, $u_{m,k}^j$ belongs to $H_\alpha(0,T;L^2(0,\ell))\cap L^2(0,T;H^{2}(0,\ell))$ and it is the unique solution of \eqref{eq1} with $b=b_j$, $q=q_j+\mu_m$, $\sigma=\sigma_j$ and $f=f_{m,k}^j$. Moreover, we have $u_{m,k}^j\in C((T-\delta,T];C^1[0,\ell])$ and \eqref{t4a} implies
\bel{t4d}\partial_s u_{m,k}^1(s,t_n)|_{s=\ell}=\partial_s u_{m,k}^2(s,t_n)|_{s=\ell},\quad n\in\mathbb N,\ m\in\mathbb N,\ k=1,\ldots,\ell_m.\ee
Since the map $g_1$ is non-uniformly vanishing, there exists $m_0\in\mathbb N$ and $k_0\in\{1,\ldots,\ell_{m_0}\}$ such that $f_{m_0,k_0}^1\not\equiv0$. Therefore, from Theorems \ref{t1} and \ref{t2}, condition \eqref{t4d}, with $m=m_0$ and $k=k_0$, implies that \eqref{t1ee} and \eqref{t2a} are fulfilled. In particular, we have $V_1=V_2$, $\sigma_2=\beta\sigma_1$  and, in view of Propositions \ref{p1} and \ref{p2} and condition \eqref{t4d}, we have
$$\partial_s \tilde{u}_{m,k}^1(s,t_n)|_{s=\ell}=\partial_s \tilde{u}_{m,k}^2(s,t_n)|_{s=\ell},\quad n\in\mathbb N,\ m\in\mathbb N,\ k=1,\ldots,\ell_m,$$
with $\tilde{u}_{m,k}^j$, $m\in\mathbb N$, $k=1,\ldots,\ell_m$  and $j=1,2$, the solution of \eqref{eq11} with $V=V_1+\mu_m$, $\sigma=\sigma_j$ and
$$\tilde{f}(s)=e^{-\frac{1}{2}\int_\ell^{s}b_j(\tau)\d\tau}f_{m,k}^j(s),\quad s\in(0,\ell).$$
Thus repeating the argument at the end of the proof of Theorem \ref{t1} yields
$$e^{-\frac{1}{2}\int_\ell^{s}b_1(\tau)\d\tau}f_{m,k}^1(s)=\beta e^{-\frac{1}{2}\int_\ell^{s}b_2(\tau)\d\tau}f_{m,k}^2(s),\quad m\in\mathbb N,\ k=1,\ldots,\ell_m,\ s\in(0,\ell).$$
Finally, by noting that, for all $m\in\mathbb N$, we have
$$e^{-\frac{1}{2}\int_\ell^{x_1}b_j(\tau)\d\tau}f_{m,k}^j(x_1)=\int_\omega \left(e^{-\frac{1}{2}\int_\ell^{x_1}b_j(\tau)\d\tau}g_j(x_1,x')\right)\psi_{m,k}(x')\d x',\quad  k=1,\ldots,\ell_m,\ s\in(0,\ell)$$
and recalling that $\{\psi_{m,k}:\ m\in\mathbb N,\ k=1,\ldots,\ell_m\}$ forms an orthonormal basis in $L^2(\omega)$, we deduce that \eqref{t4b} is also fulfilled. This completes the proof of the theorem. \end{proof}

\section{Numerical experiments and discussions}\label{sec:numer}
In this section, we present numerical results to show the feasibility of recovering multiple parameters, including the convection coefficient $b$, the potential $q$, the source $f$, and the source strength $\sigma$, in the subdiffusion model from the flux observation $g$ at the end point $x=1$. We employ the Levenberg-Marquardt (LM) method \cite{Levenberg:1944,Marquardt:1963}, which iteratively refines the parameter estimates. We describe the method for recovering $b$ and $f$, and the other cases are similar. Let $F: (b, f) \mapsto \partial_x u|_{x=1} $  denote the forward operator. At the \(k\)th iteration, we linearize the operator \(F\) at the current estimate \((b_k, f_k)\), and construct the following regularized least-squares problem:
\[
J_k(b, f) = \left\| F(b_k, f_k) - g + \nabla_b F \cdot (b - b_k) + \nabla_f F \cdot (f - f_k) \right\|_2^2 + \frac{\gamma_{b,k}}{2} \|b - b_k\|_{H^1}^2 + \frac{\gamma_{f,k}}{2} \|f - f_k\|_{H^1}^2,
\]
where the parameters $\gamma_{b,k}$ and $\gamma_{f,k}$ are decreased geometrically during iterations.
The updated tuple \( (b_{k+1}, f_{k+1}) \) is then given by
\[
(b_{k+1}, f_{k+1}) = \arg\min_{b,f} J_k(b, f).
\]

First consider the one-dimensional fractional PDE on the domain $(0,1)$ with
$T=0.1$:
\[
\partial^\alpha_t u - \partial^2_x u + b(x)\partial_x u + q(x)u = \sigma(t)f(x),
\]
equipped with the boundary condition \(u(0,t) = u(1,t) = 0\), and initial condition \(u(x,0) = 0\). The exact coefficients $b^\dag$, $q^\dag$, $f^\dag$ and $\sigma^\dag$ are given by
\begin{align*}
b^\dag(x) &=
\begin{cases}
16x^3-12x^2+3x & \text{if } x < 0.25, \\
0.25 & \text{if } x \geq 0.25,
\end{cases}, \quad
q^\dag(x) = e^{-x},\\
f^\dag(x) &= x\sin(2\pi x) \cdot \mathbf{1}_{[0,0.5)}(x),
\quad \sigma^\dag(t) = (50 - 1000t) \cdot \mathbf{1}_{[0,0.05]}(t).
\end{align*}
Then consider the following three cases, with the measurement being the Neumann boundary data at the point \(x = 1\).
\begin{itemize}
    \item[(a)] \(\alpha = 0.5\),  \(q\) and \(\sigma\) are known over $[0,1]$, and \(b\) is known over  \([0.25, 1]\). Recover \(b\) over \([0, 0.25)\) and \(f(x)\) over \([0, 0.5]\).
    \item[(b)] $\alpha=0.5$, \(b\) and \(\sigma\) are known over \([0 , 1]\), and \(q\) is known over \([0.25, 1]\). Recover \(q\) over \([0, 0.25]\) and the source \(f(x)\) supported over \([0, 0.5]\).
    \item[(c)] \(\alpha = 3/(2\pi)\), \( q \) is known over $[0,1]$, and \( b \) over \([0.25, 1]\). Recover \( b \) over \([0, 0.25]\),  \( f \) supported over \([0, 0.5]\), and \(\sigma\) supported on \([0, 0.05]\).
\end{itemize}

Cases (a) and (b) are to validate the two parts of Theorem \ref{t1_intro}, and case (c) is to validate  Theorem \ref{t2_intro}. We generate the exact data using \(2,\!048\) spatial / temporal nodes, and then solve the inverse problem using  \(1,\!024\) spatial and \(2,\!048\) temporal nodes. To measure the accuracy of the reconstructions, we use the relative errors, e.g., $e_b=\|b^\dag-\hat b\|_2/\|b^\dag\|_2$ for $\hat b$.

\begin{table}[htb!]
\centering\setlength{\tabcolsep}{3pt}
    \caption{The numerical results for different cases at various noise levels.}\label{tab:results}
    \begin{tabular}{ccc|cccc}
    \toprule
    \multicolumn{3}{c}{(a) }&\multicolumn{4}{c}{(c)}\\
    \cmidrule(lr){1-3} \cmidrule(lr){4-7}
        \(\epsilon\) & $e_b$ & $e_f$ & $\epsilon$ &\(e_b\) & \(e_f\) & \(e_\sigma\)\\
        \midrule
        \(0\) & 0.48\% & 2.24\% & $0$ &0.46\% & 2.13\% & \(5.60 \times 10^{-9}\)  \\
        \(10^{-5}\) & 0.53\% & 4.63\% & $10^{-5}$ & 0.47\% & 4.21\% & \(1.13 \times 10^{-6}\)\\
        \(10^{-4}\) & 1.25\% & 4.76\% & $10^{-4}$ & 1.06\% & 5.03\% & \(4.61 \times 10^{-5}\)\\
    \midrule
    \multicolumn{3}{c}{(b)}&\multicolumn{4}{c}{(d)}\\
    \cmidrule(lr){1-3} \cmidrule(lr){4-7}
            \(\epsilon\) & \(e_q\) & \(e_f\) & \(\epsilon\) & \(e_b\) & \(e_f\) & \(e_\sigma\) \\
        \midrule
        \(0\) & 0.23\% & 3.60\% & \(0\) & 0.68\% & 1.90\% & \(2.61 \times 10^{-4}\)\\
        \(10^{-5}\) & 0.27\% & 3.90\% & \(10^{-4}\) & 1.68\% & 1.21\% & \(1.73 \times 10^{-4}\)\\
        \(10^{-4}\) & 3.77\% & 9.15\% & \(10^{-3}\) & 49.2\% & 7.79\% & \(7.60 \times 10^{-3}\)\\
    \bottomrule
    \end{tabular}

\end{table}

\begin{figure}[hbt!]
    \centering
    \begin{tabular}{cc}
         \includegraphics[width = 0.4\linewidth]{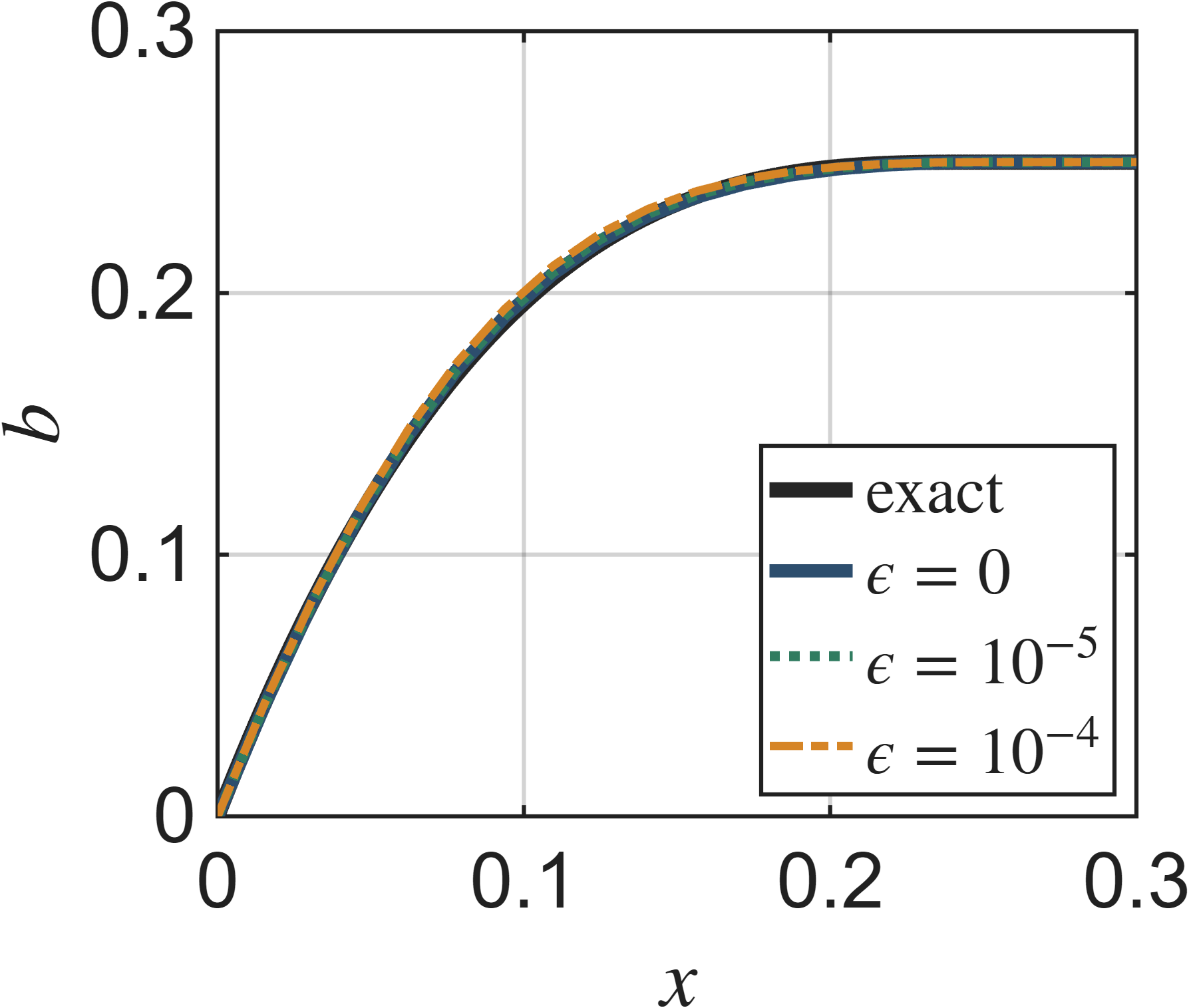}&  \includegraphics[width = 0.4\linewidth]{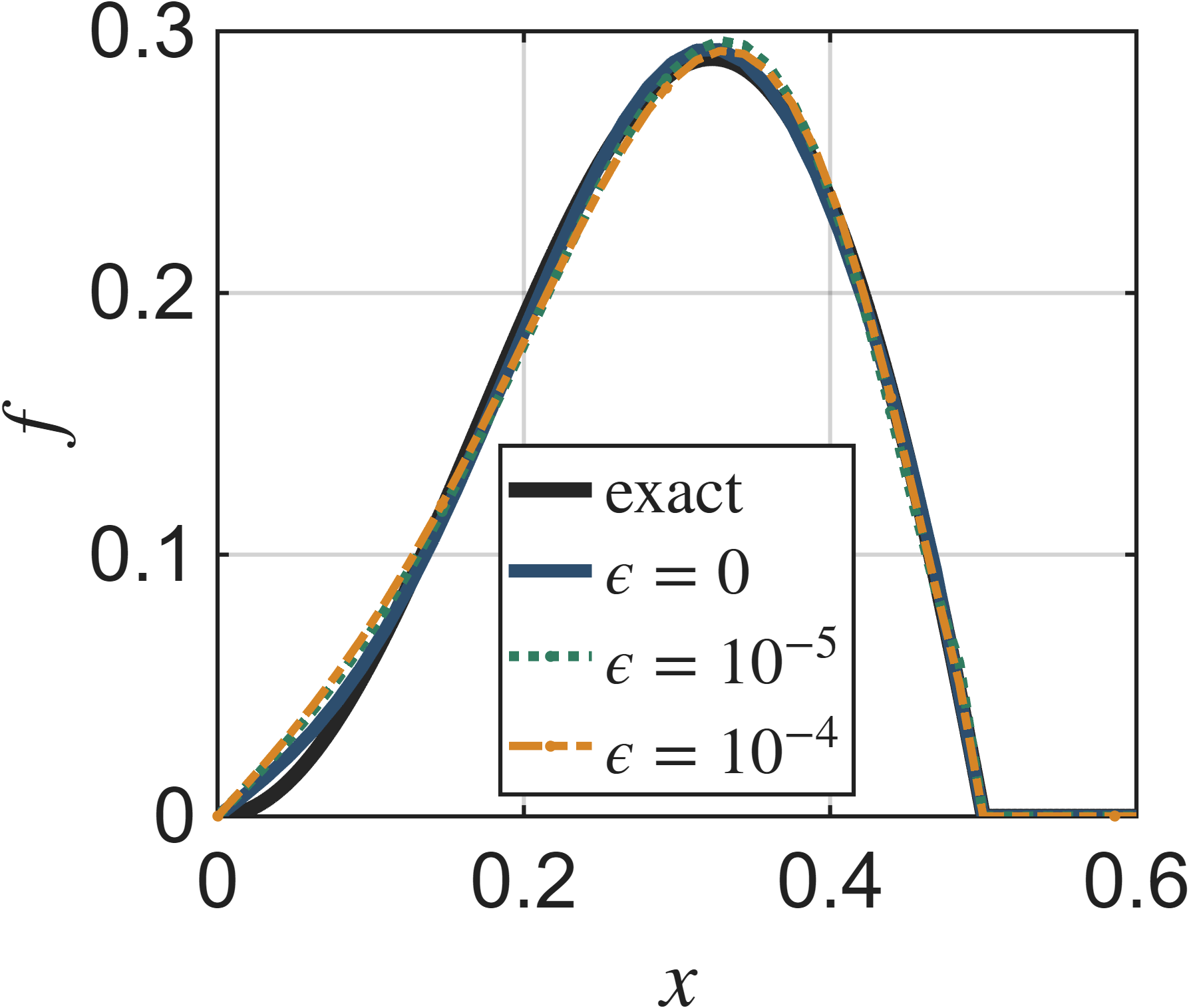}\\
         (a) $b$ & (b) $f$
    \end{tabular}
    \caption{The recovered \(b\) and \(f\) for case (a) at different noise levels.  }
    \label{fig:coefficients}
\end{figure}

\begin{figure}[hbt!]
    \centering \setlength{\tabcolsep}{0pt}
     \begin{tabular}{ccc}
  \includegraphics[width=0.33\linewidth]{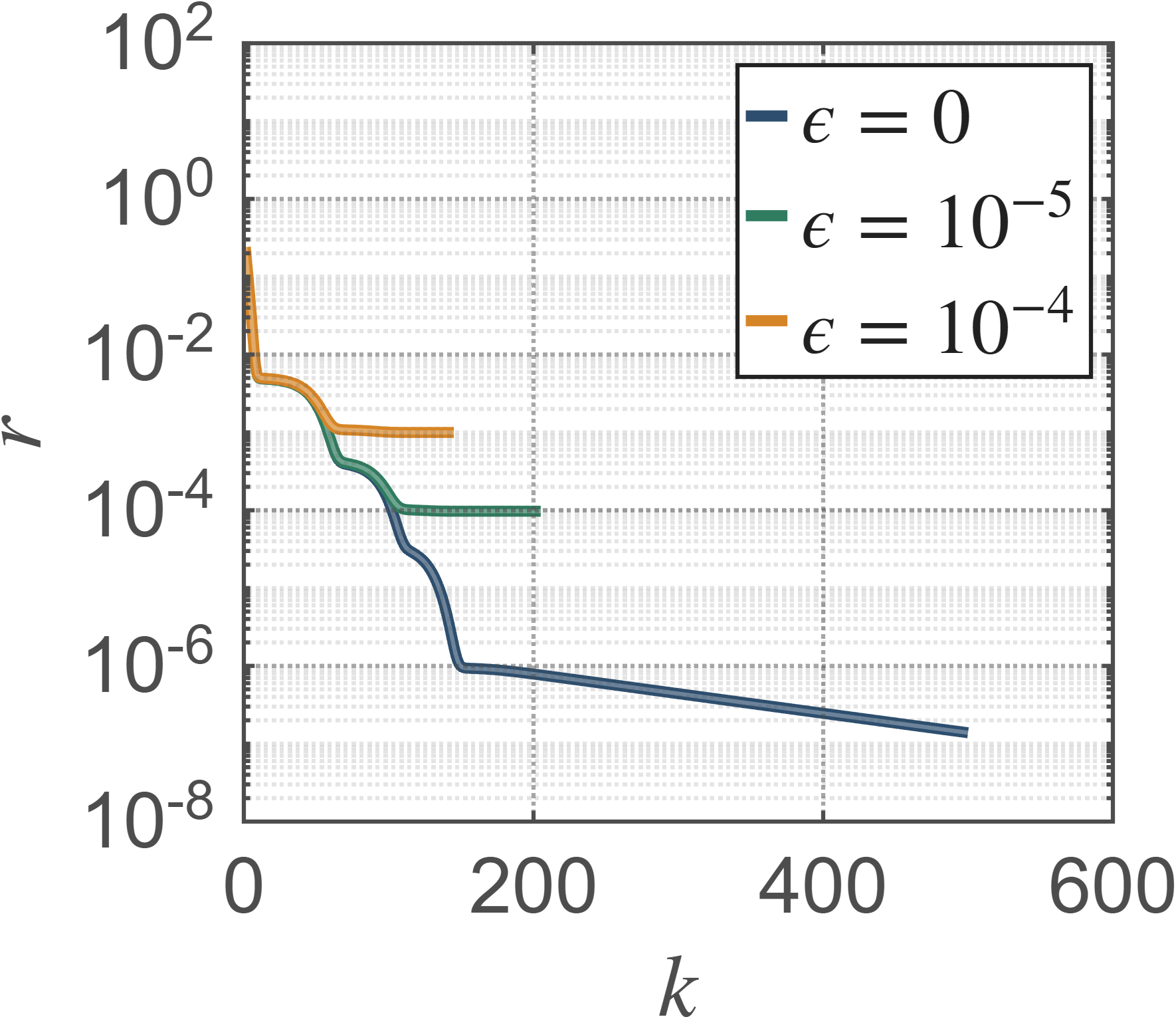}&\includegraphics[width=0.33\linewidth]{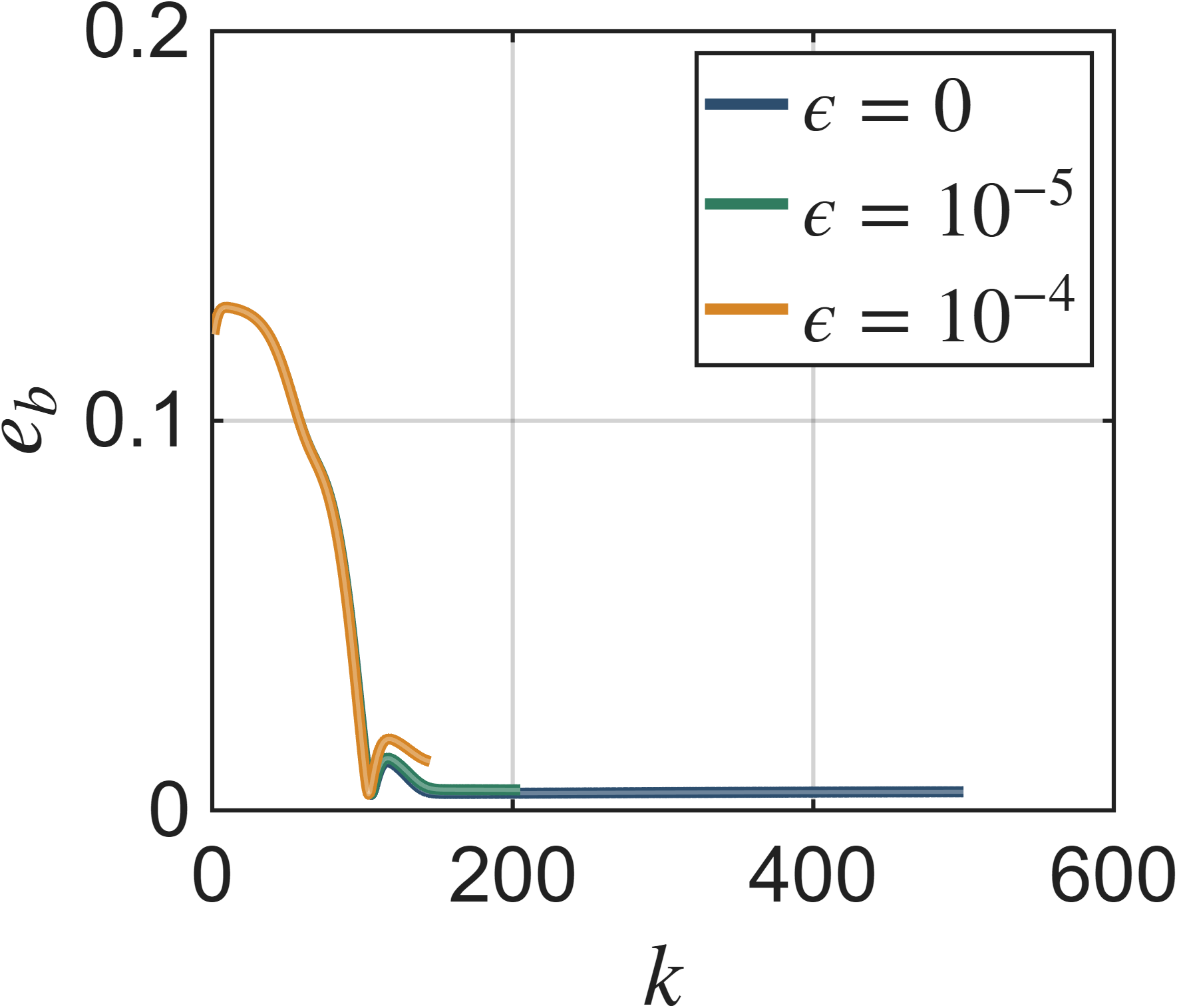}&\includegraphics[width=0.33\linewidth]{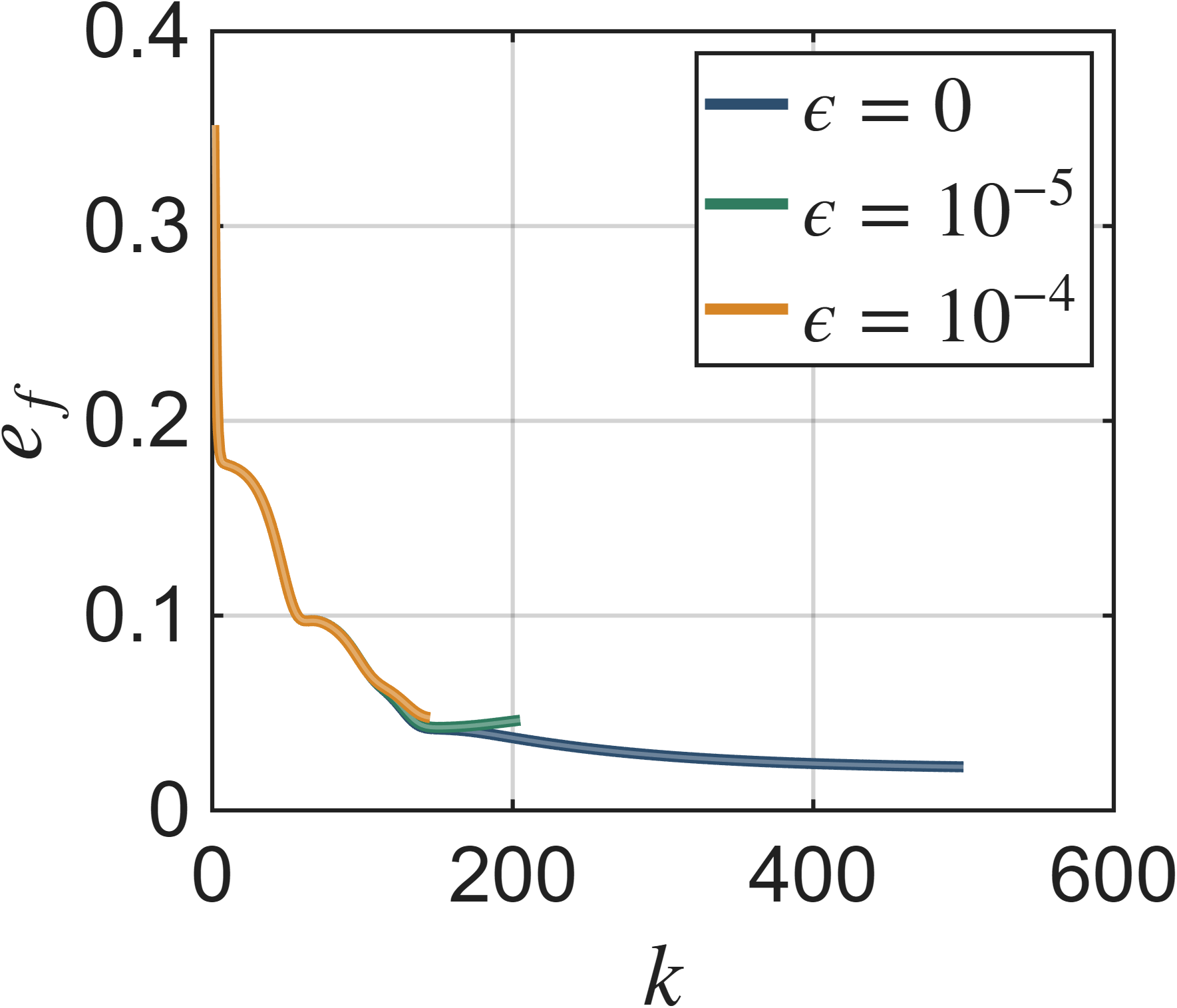}\\
  (a) $r$  & (b) $e_b$ & (c) $e_f$
  \end{tabular}
    \caption{The convergence of the Levernberg-Marquardt method for case (a).  }
    \label{fig:loss}
\end{figure}

In case (a), to enforce the \(C^1\) continuity of \(b\) at $x=0.25$, we parameterize \(b = b_0(x) + \beta(x)(x - 0.25)\) and optimize \(\beta\), with an initial guess \(b_0 \in C^1\). Fig. \ref{fig:coefficients} shows the recovered \(b\) and \(f\) at different noise levels, and indicates that the reconstructions \(\hat b\) and $\hat f$ are accurate for up to \(\epsilon = 10^{-4}\). Moreover, the relative error for \(b\) is below $1\%$ for exact data, clearly showing the feasibility of the numerical reconstruction. Fig. \ref{fig:loss} shows the convergence of the residual $r= \|F(b_k, f_k) - g\|_2$ and relative errors, showing the efficiency of the LM algorithm. Table \ref{tab:results}(a) shows that \(f\) is less sensitive to the noise than \(b\) (e.g., \(e_f\) increases by 2\% at \(\epsilon = 10^{-4}\) vs. 136\% for \(e_b\)), since the source $f$ typically has stronger influence on the boundary fluxes

\begin{figure}[hbt!]
    \centering
    \begin{tabular}{cc}
        \includegraphics[width=0.4\linewidth]{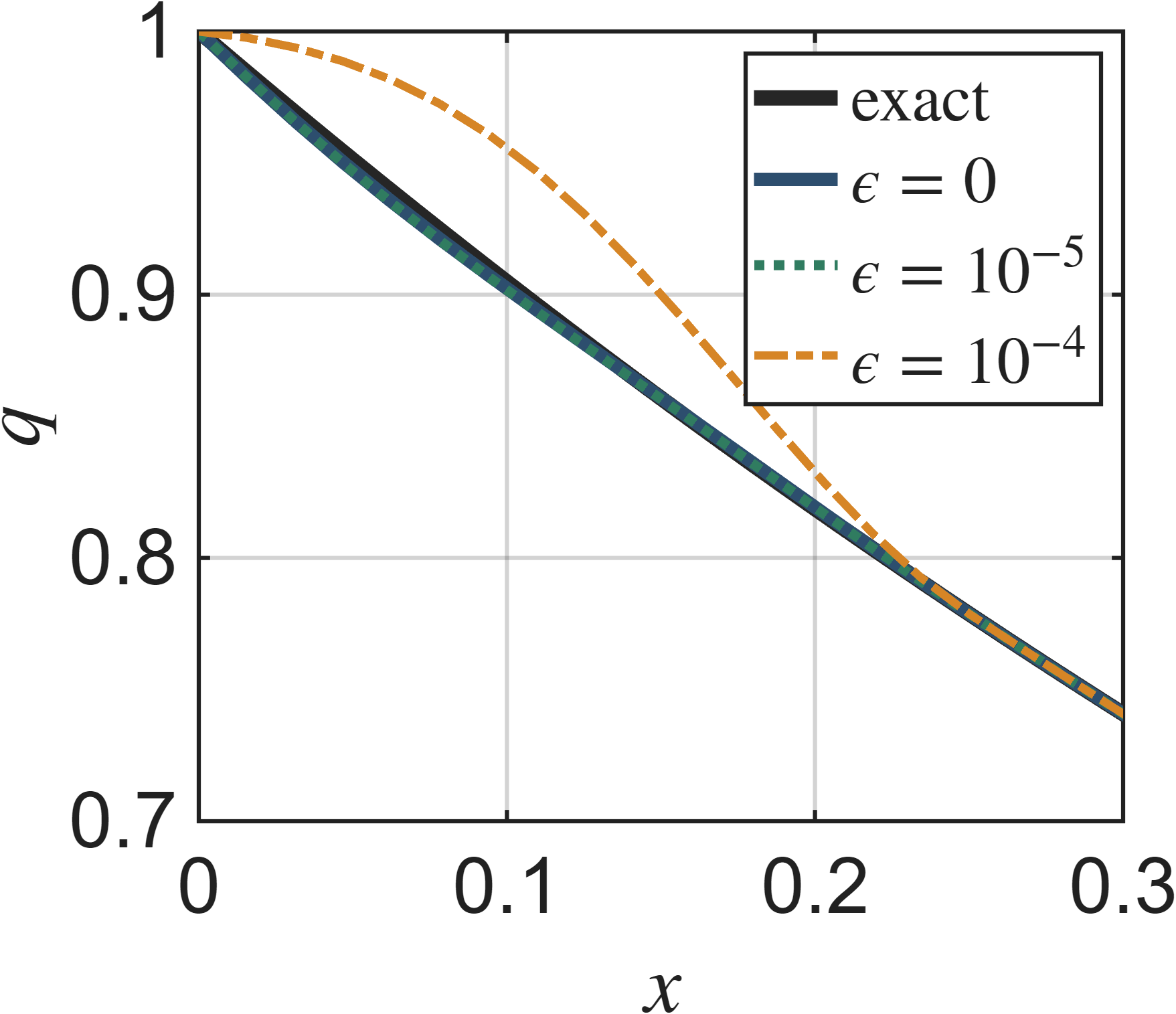} &
        \includegraphics[width=0.4\linewidth]{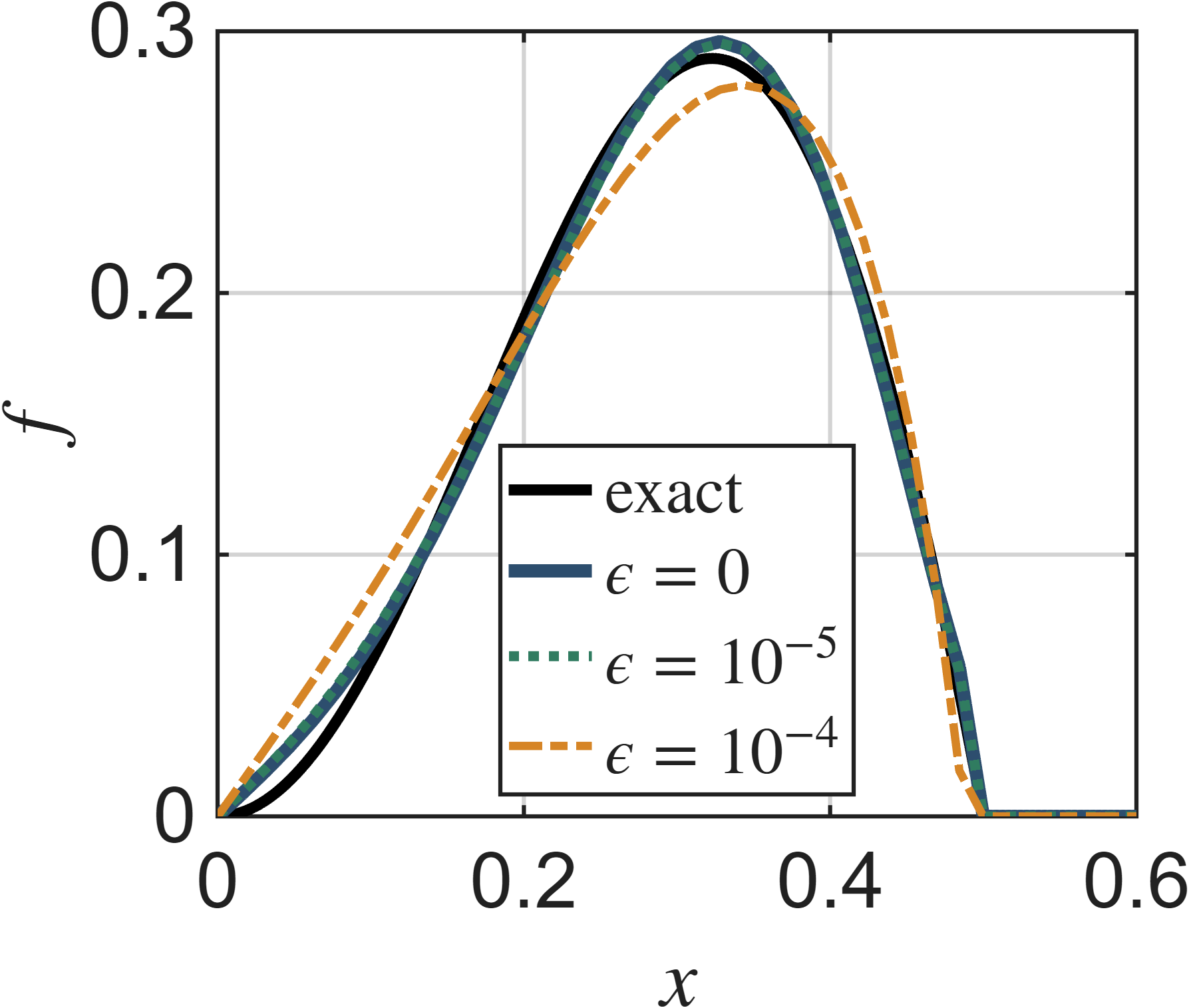}\\
        (a) $q$ & (b) $f$
    \end{tabular}
    \caption{The recovered \(q\) and \(f\) for case (b) at three noise levels.}
    \label{fig:coefficients_q}
\end{figure}

\begin{figure}[hbt!]
    \centering\setlength{\tabcolsep}{0pt}
    \begin{tabular}{ccc}
        \includegraphics[width=0.33\linewidth]{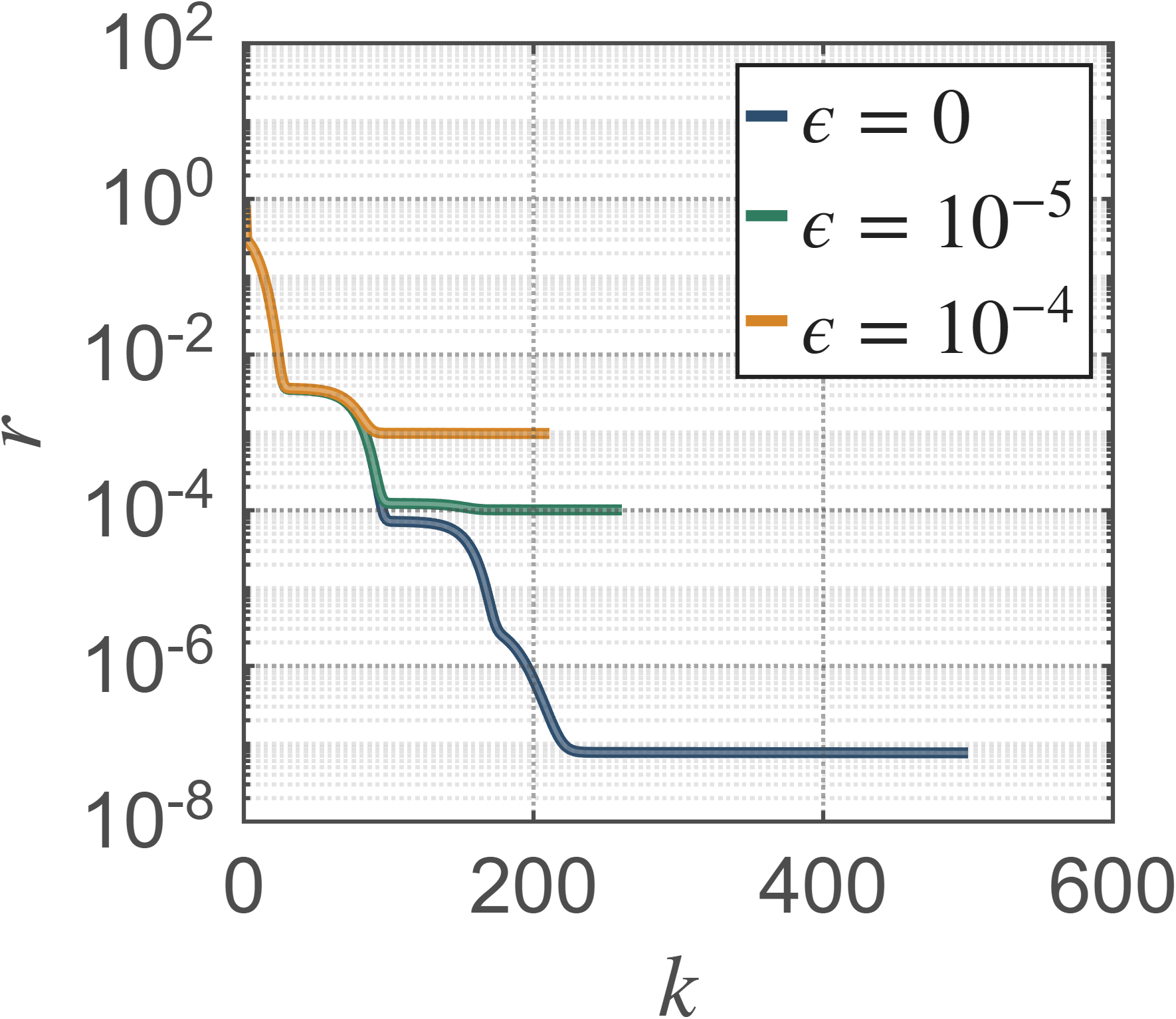} &
        \includegraphics[width=0.33\linewidth]{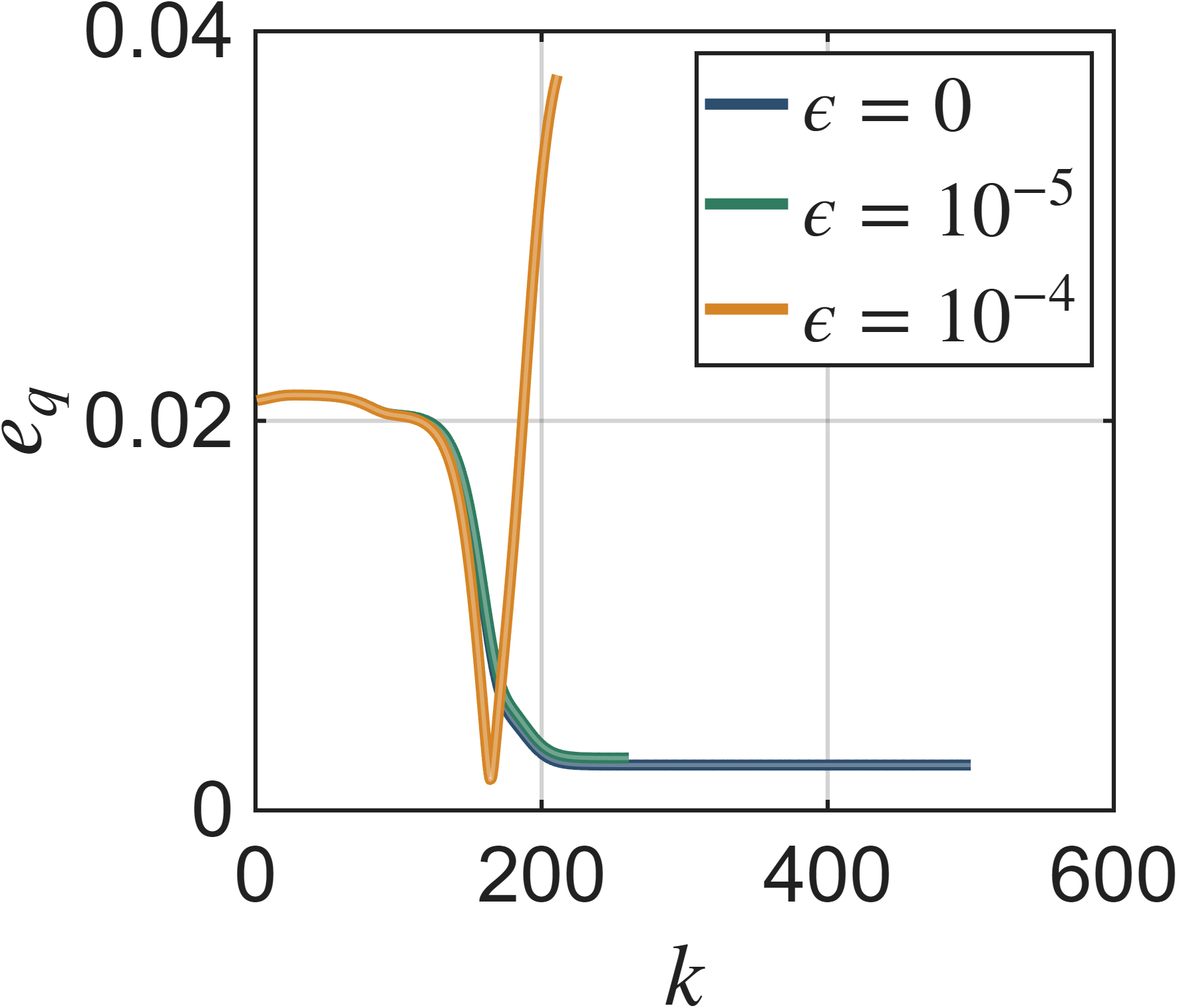} &
        \includegraphics[width=0.33\linewidth]{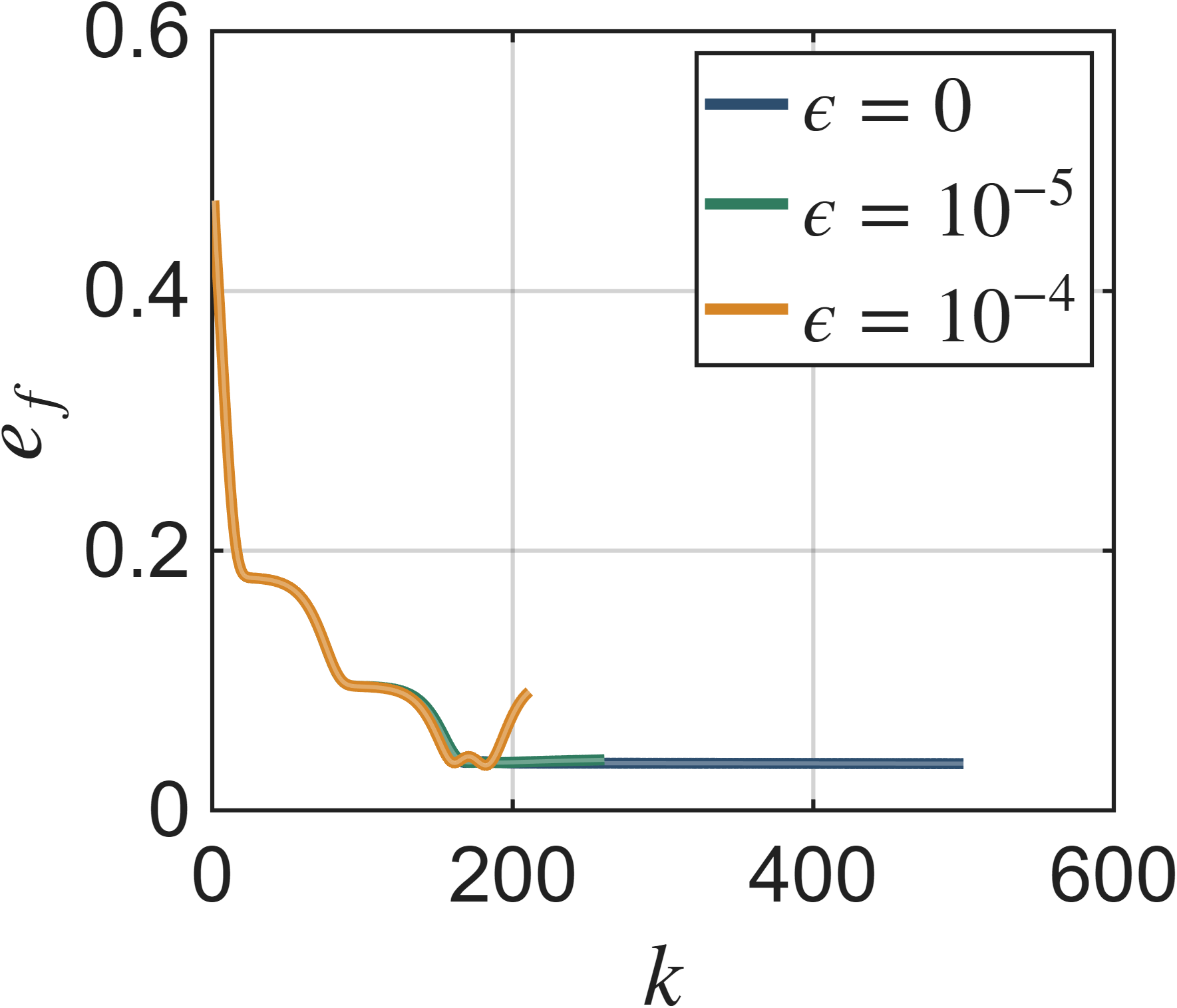}\\
        (a) $r$ & (b) $e_q$ & (c) $e_f$
    \end{tabular}
    \caption{The convergence of the Levenberg-Marquadt method for case (b).}
    \label{fig:loss_q}
\end{figure}

The results for case (b) are given in Table \ref{tab:results}(b), and Figs. \ref{fig:coefficients_q} and \ref{fig:loss_q}. Compared with case (a), Table \ref{tab:results}(b) shows that \(q\) can be recovered more robustly than \(b\) (e.g., \(e_q = 0.23\%\) versus \(e_b = 0.48\%\) for exact data). At \(\epsilon = 10^{-4}\), \(e_q\) increases sharply (15×) while \(e_f\) grows moderately (2.5×), suggesting higher susceptibility of $q$ to the noise. Fig. \ref{fig:coefficients_q} shows that the recovered \(q\) are nearly identical with the exact one up to the case \(\epsilon = 10^{-5}\), and the recovered \(f\) can capture the location and size of the bump at all noise levels.

\begin{figure}[hbt!]
    \centering \setlength{\tabcolsep}{0pt}
    \begin{tabular}{ccc}
        \includegraphics[width=0.33\linewidth]{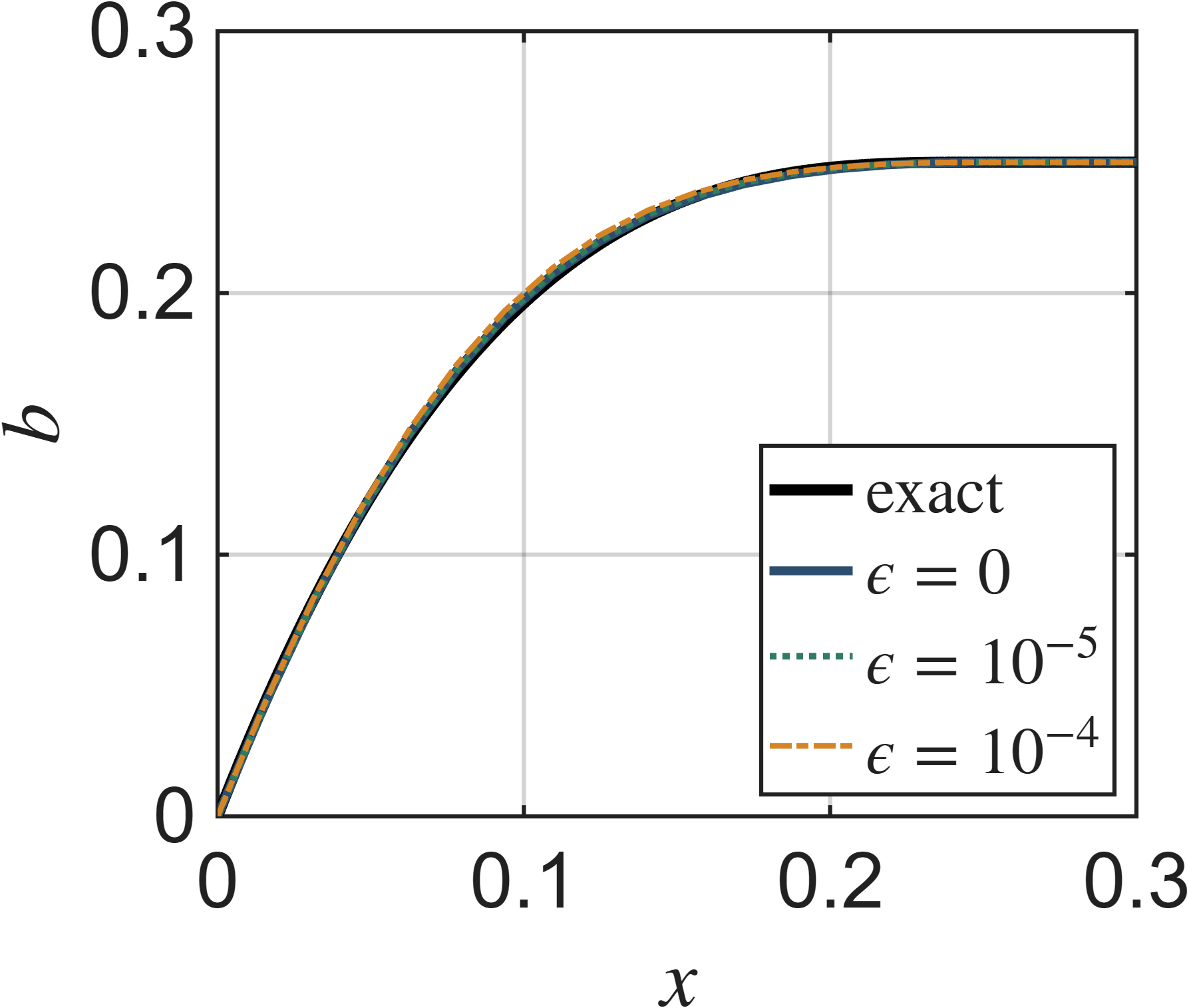} &
        \includegraphics[width=0.33\linewidth]{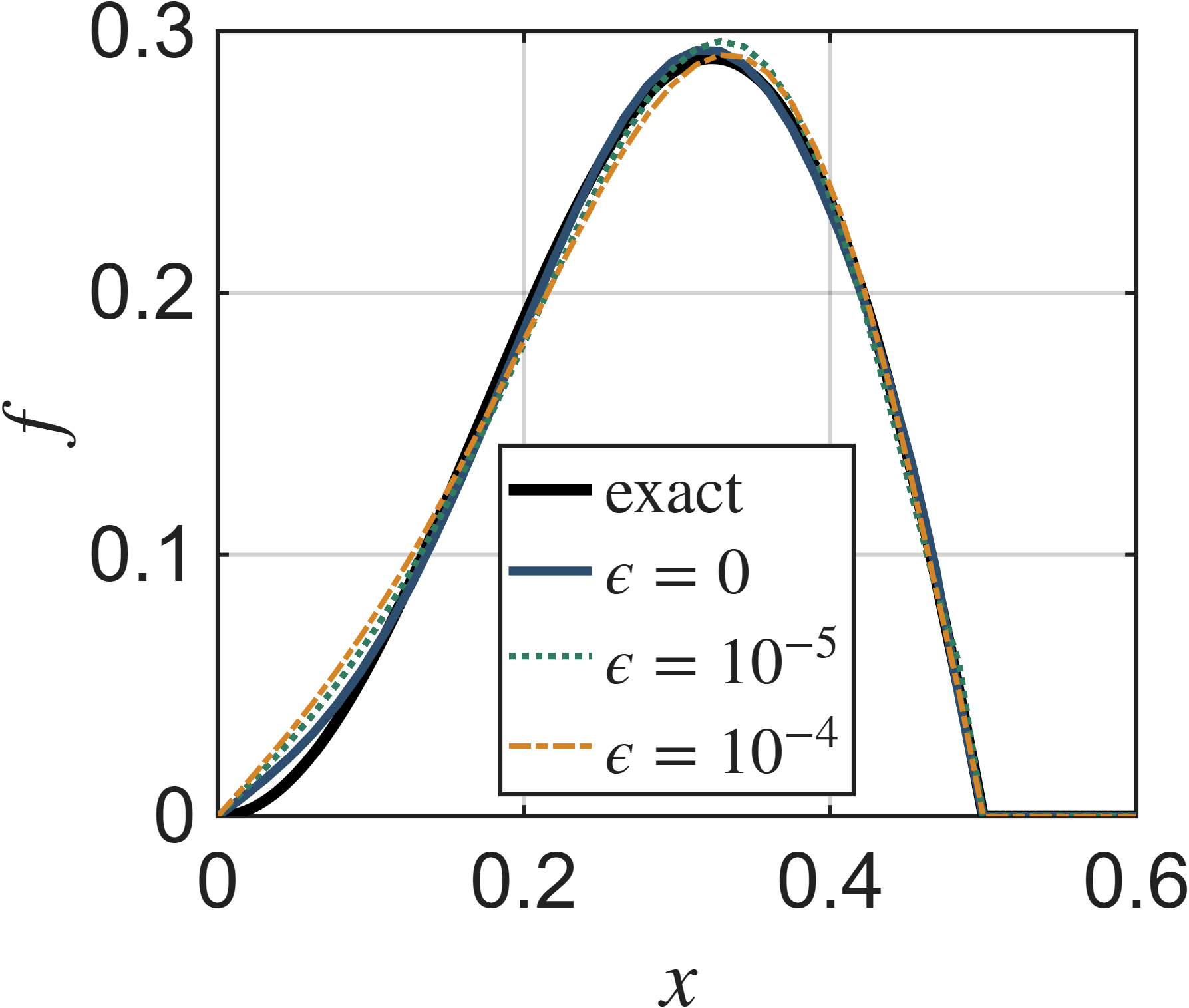} &
        \includegraphics[width=0.33\linewidth]{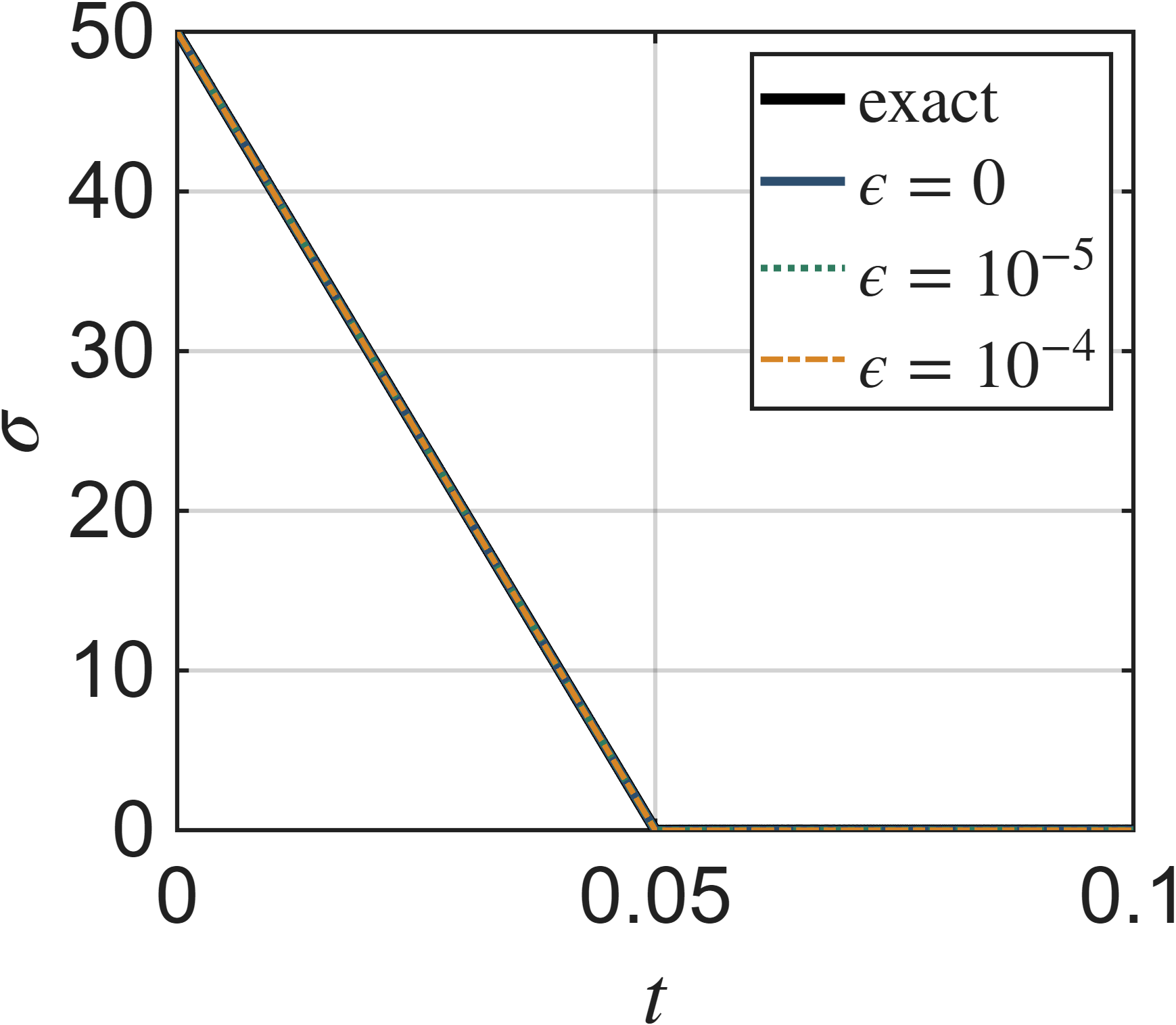}\\
        (a) $b$ & (b) $f$ & (c) $\sigma$
    \end{tabular}
    \caption{The recovered \(b\), \(f\), and \(\sigma\) for case (c) at three noise levels.  }
    \label{fig:coefficients_tri}
\end{figure}

The results for case (c) are presented in Fig. \ref{fig:coefficients_tri} and Table \ref{tab:results}(c). Note that \(\sigma(t)\) is recovered with very high accuracy (\(e_\sigma < 10^{-4}\) for \(\epsilon = 10^{-4}\)), since the source component  $\sigma $ is well aligned with the flux measurement. The relative errors \(e_b\) and \(e_f\) are largely comparable with that in Tables \ref{tab:results}(a) and (b). The recovered \(\sigma\)  are visually distinguishable from the exact one, and the recovered \(f\)  and $b$ are comparable with the preceding experiments, cf. Fig. \ref{fig:coefficients_tri}.

\begin{figure}[hbt!]
    \centering
    \begin{tabular}{cccc}
        \includegraphics[width=0.4\linewidth]{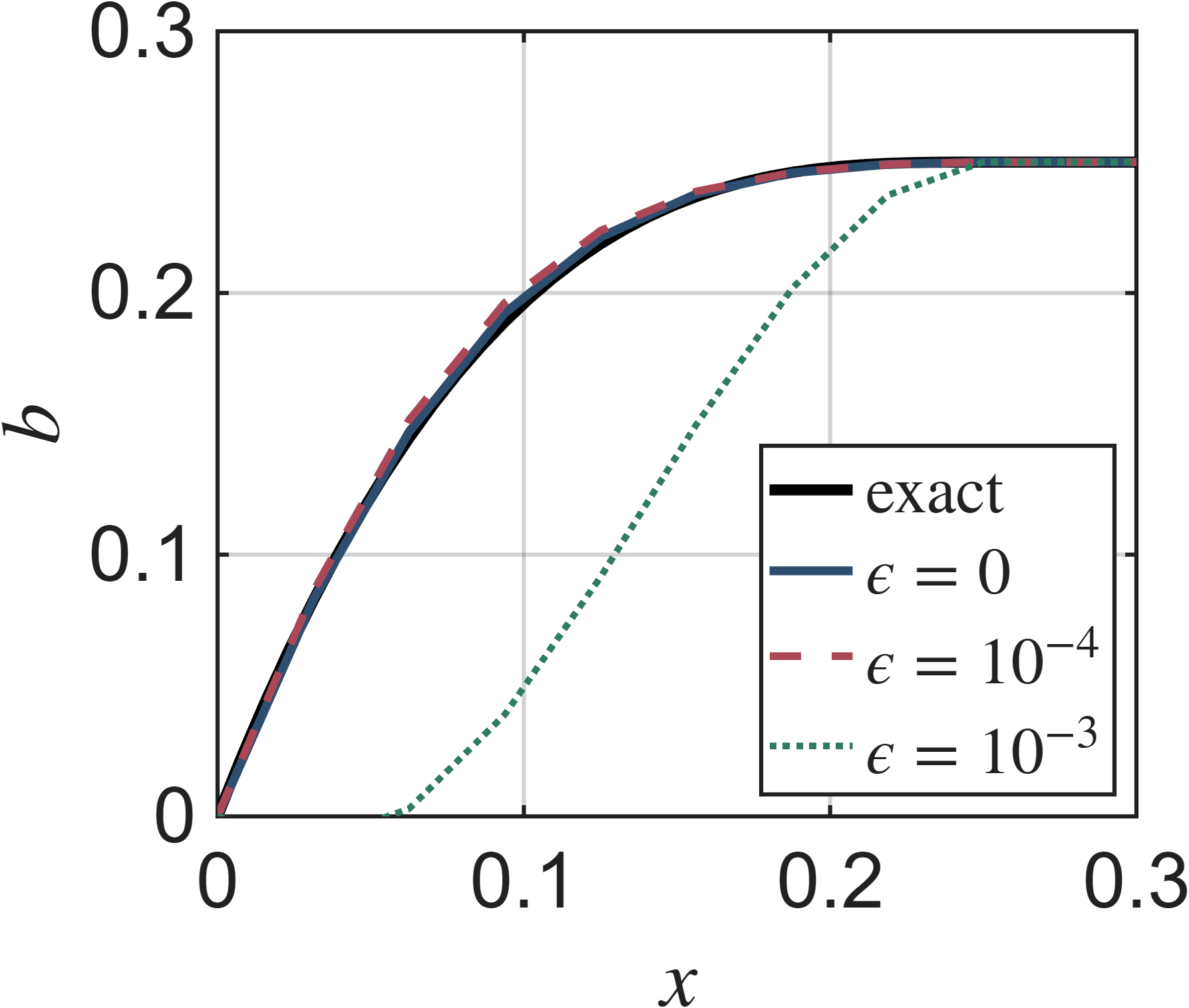} &
        \includegraphics[width=0.4\linewidth]{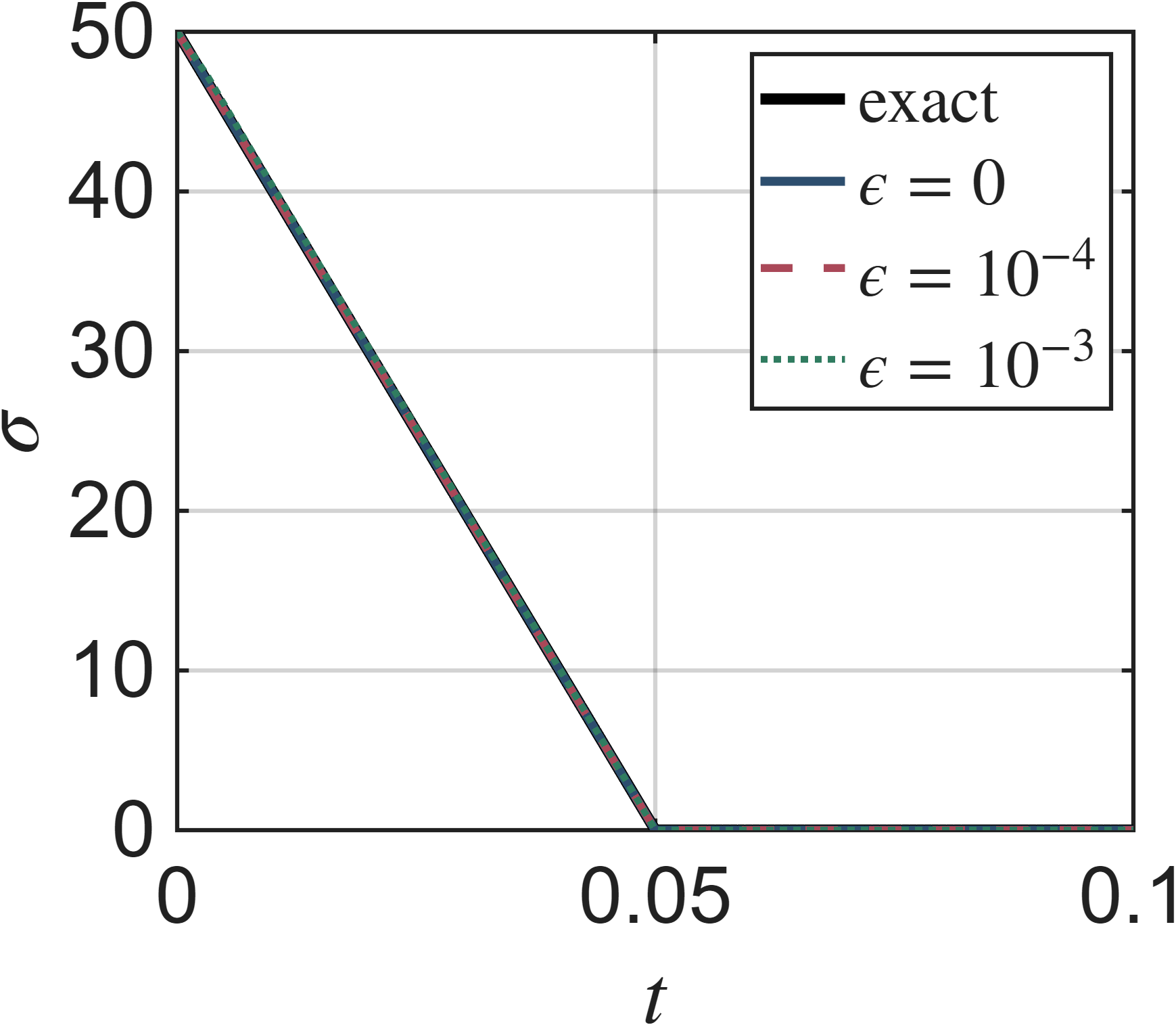}  \\
        (a) $b(x_1)$ & (b) $\sigma$
    \end{tabular}
    \caption{The recovered \(b(x_1)\)  and \(\sigma\) for case (d). }
    \label{fig:bs2d}
\end{figure}

\begin{figure}[hbt!]
    \centering\setlength{\tabcolsep}{0pt}
    \begin{tabular}{cccc}
        \includegraphics[width=0.25\linewidth]{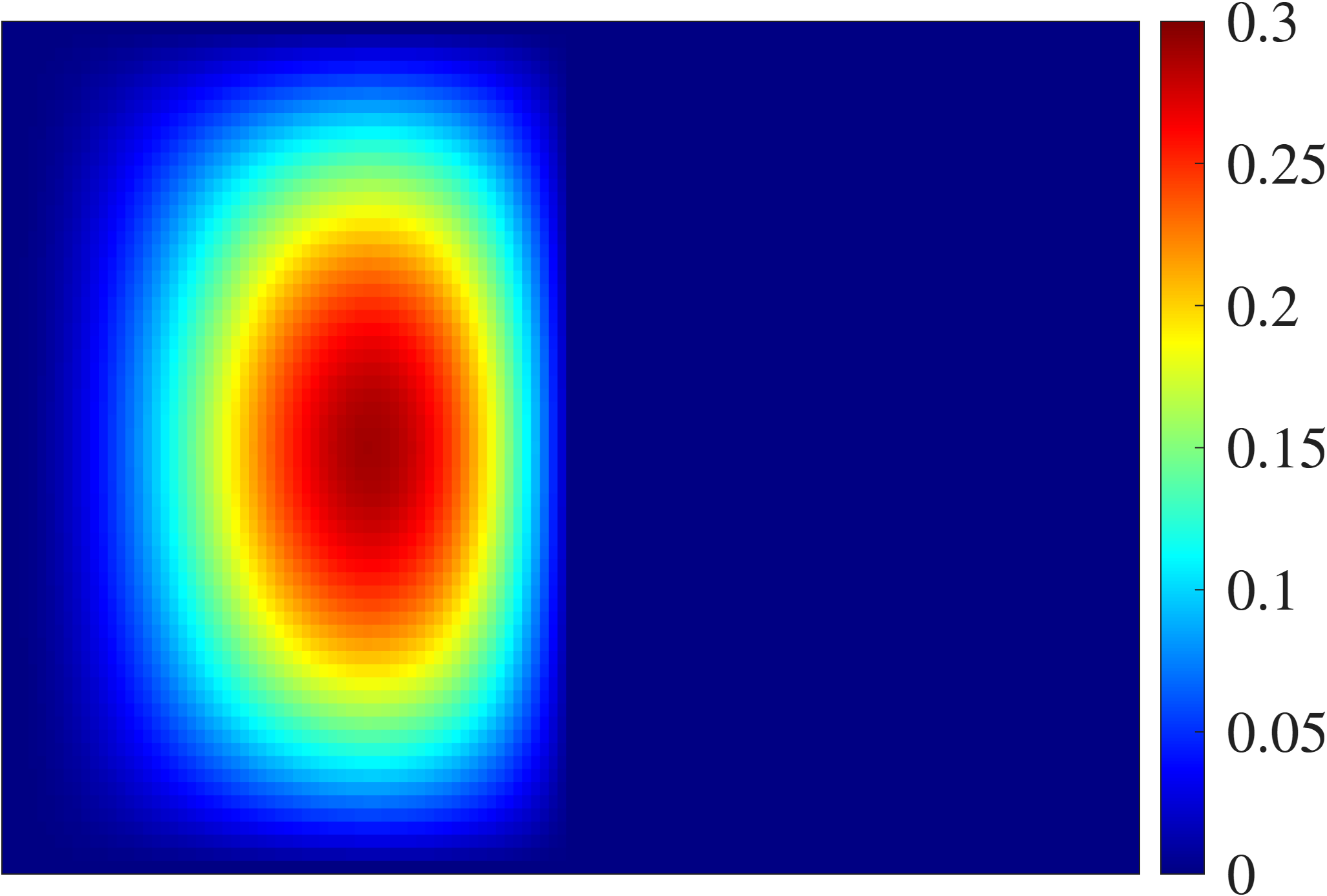} &         \includegraphics[width=0.25\linewidth]{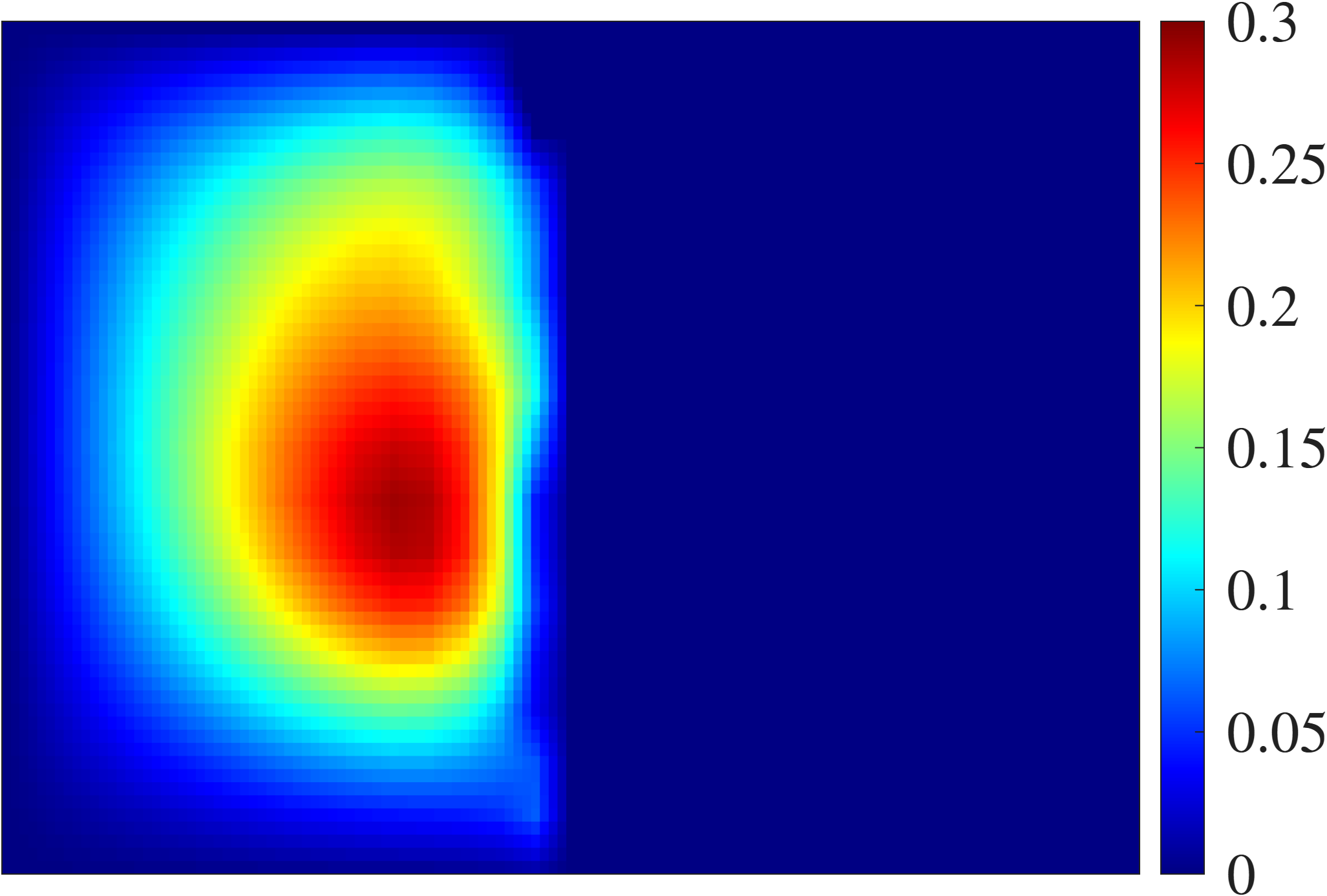} & \includegraphics[width=0.25\linewidth]{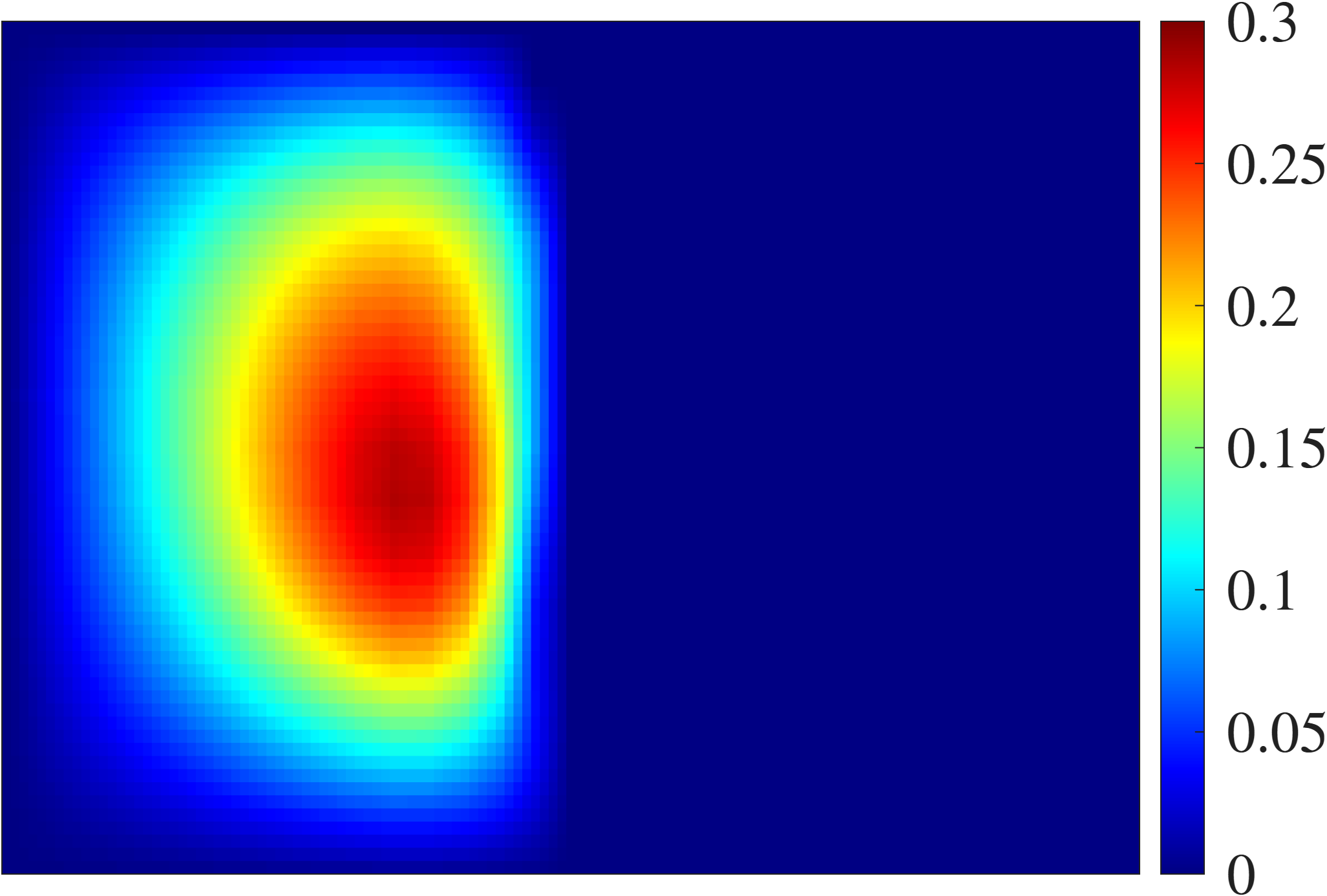} & \includegraphics[width=0.25\linewidth]{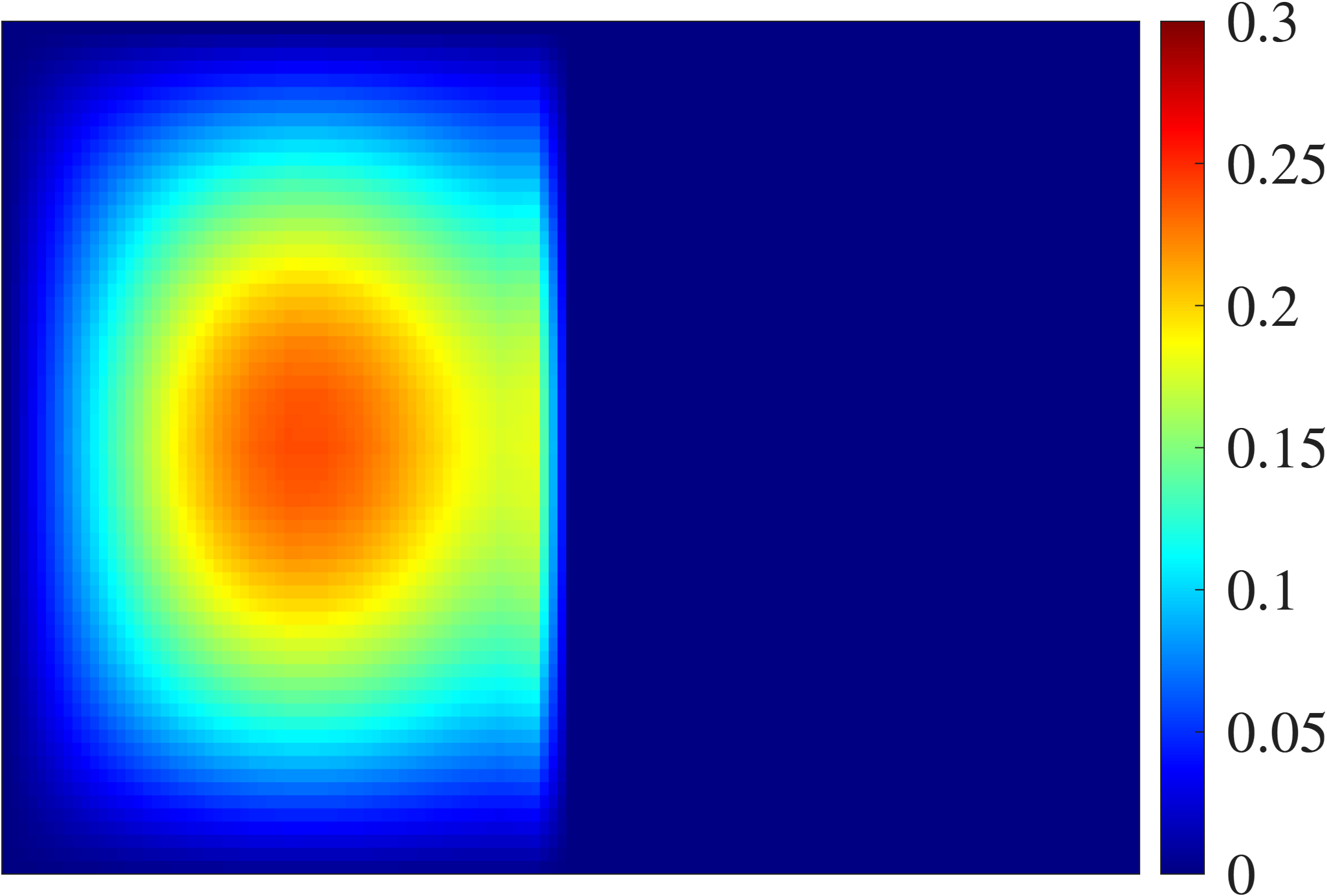} \\
&\includegraphics[width=0.25\linewidth]{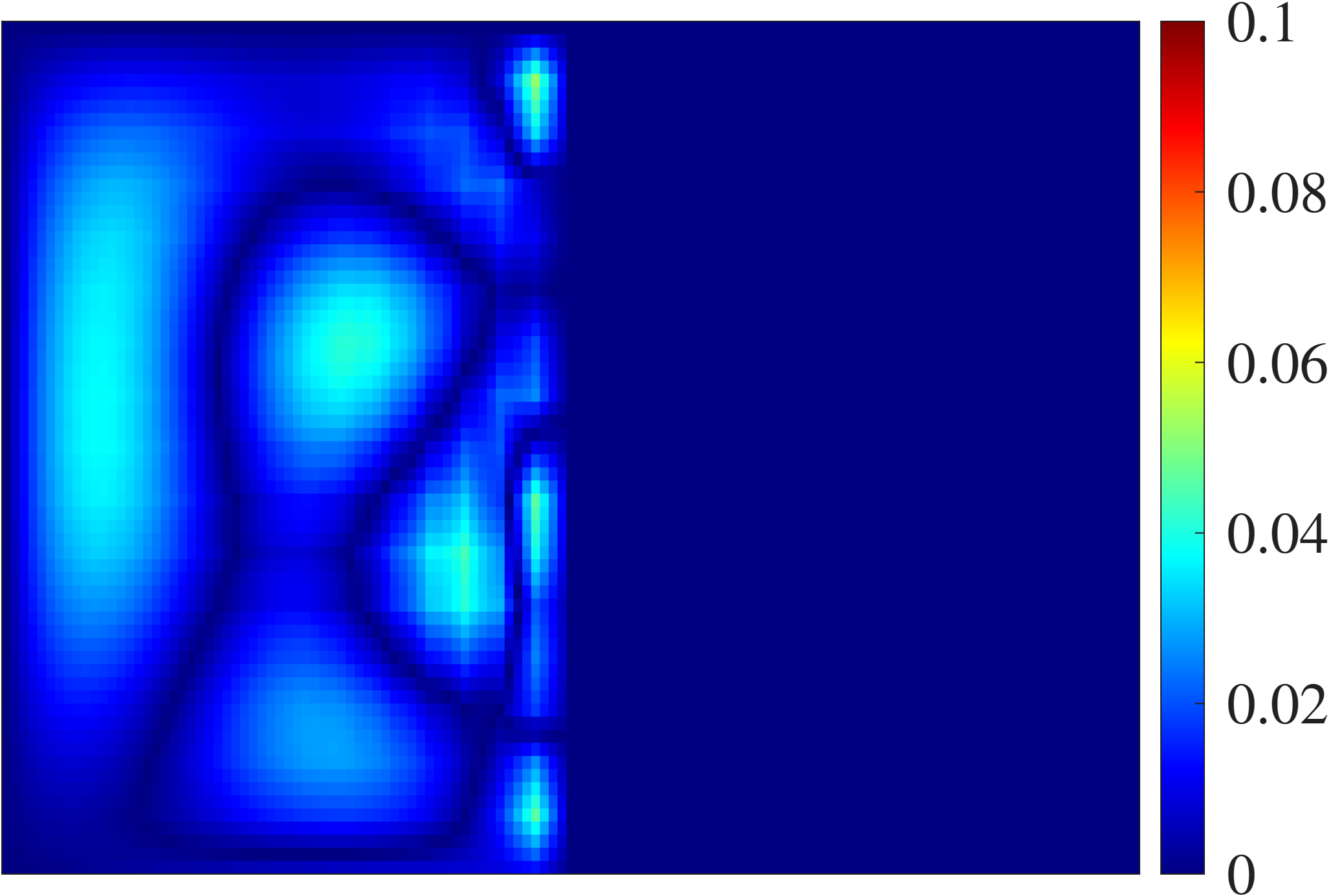} &
\includegraphics[width=0.25\linewidth]{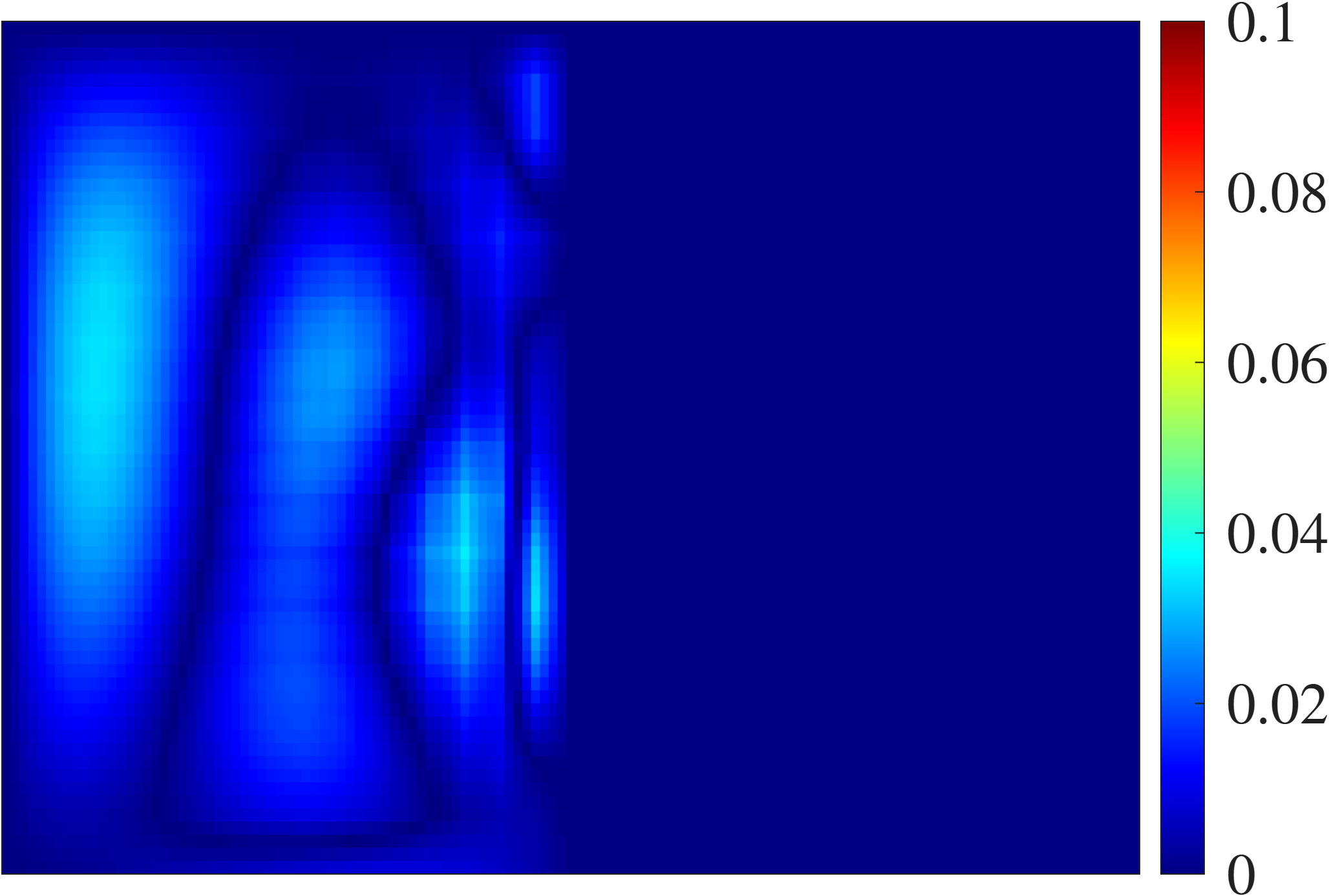} &\includegraphics[width=0.25\linewidth]{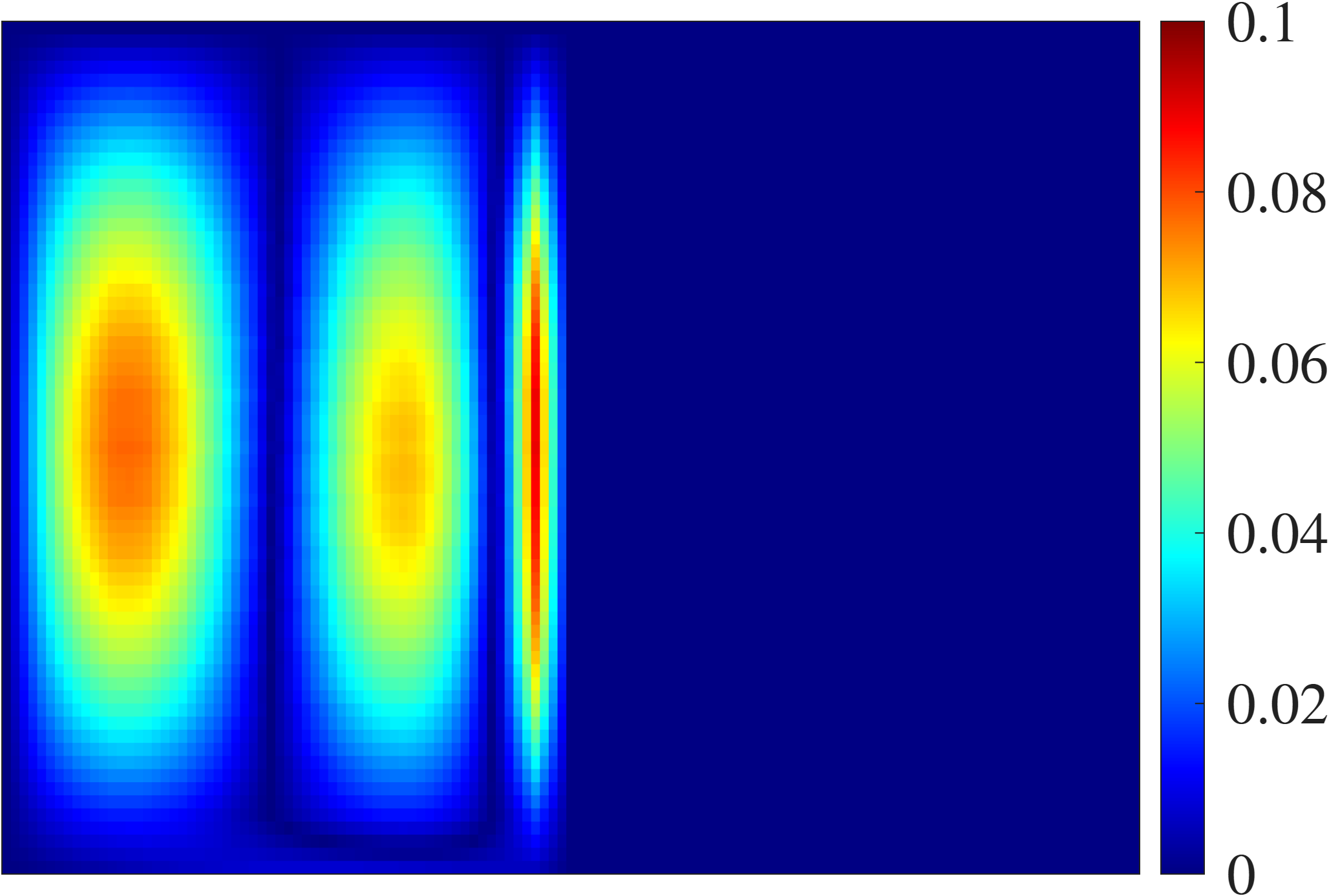}\\
(a) $f^\dag$ & (b) $\epsilon=0$ & (c) $\epsilon=10^{-4}$ & (d) $\epsilon=10^{-3}$
    \end{tabular}
    \caption{The recovered $f$ (top) and pointwise error (bottom) for case (d) at different noise levels. }
    \label{fig:f2d}
\end{figure}

Last we illustrate Theorem \ref{t4_intro} in the 2D case with the domain $\Omega=(0,1)^2$:
\begin{itemize}
    \item[(d)] \(\alpha = 3/(2\pi)\),  \( q \)  is known in the interval $[0,1]$ and \( b \) in \([0.25, 1]\). Recover \( b \) in \([0, 0.25]\),  \( f \) supported over \([0, 0.5]\times[0,1]\), as well as \(\sigma\) with its support in \([0, 0.05]\).
\end{itemize}
The results for case (d) at three noise levels are given in Table \ref{tab:results}(d) and Figs. \ref{fig:bs2d} and \ref{fig:f2d}. At \(\epsilon = 10^{-3}\), the reconstruction largely fails, in view of the relative error \(e_b = 49.2\%\), due to the severe ill-posedness of the problem of recovering the convection coefficient $b$. This is also clearly observed for $b(x_1)$ in Fig. \ref{fig:bs2d}(a). In contrast, the time component \(\sigma(t)\) can still be accurately recovered even at \(\epsilon = 10^{-3}\), with a relative error $e_\sigma = 0.76\%$.
Fig. \ref{fig:f2d} shows the recovered $f$. 

\appendix
\section{Appendix}
In this appendix, we recall some inverse spectral results for Sturm-Liouville type operators. All these results were initially stated for operators acting on $L^2(0,1)$. We give their extension to operators acting on $L^2(0,\ell)$. Throughout, for $j=1,2$, we fix $V_j\in L^\infty(0,\ell)$ and $A_j=-\partial_x^2 u+V_j$ acting on the space $L^2(0,\ell)$ with its domain $D(A_j)=H_0^1(0,\ell)\cap H^2(0,\ell)$, $(\lambda_{j,k})_{k\in\mathbb N}$ the increasing sequence of simple eigenvalues of $A_j$ and we associate with these eigenvalues a corresponding orthonormal basis $(\phi_{j,k})_{k\in\mathbb N}$ of eigenfunctions of $L^2(0,\ell)$.
Following \cite[Lemma 3.3]{Feizmohammadi2025}, we have the following lemma.
\begin{lem}\label{l1} We have
\bel{l1a}\lambda_{j,k}=\left(\frac{k\pi}{\ell}\right)^2+\frac{1}{\ell} \int_0^\ell V_j(y)\d y+\underset{k\to+\infty}{ o}(1) ,\ee
\bel{l1b}\limsup_{\mu\to+\infty} \frac{m(\mu)}{\sqrt{\mu}}\leq \frac{r}{\pi}.\ee
\end{lem}
\begin{proof}
Note that for $\phi\in H^2(0,\ell)$, $V\in L^\infty(0,\ell)$ and $\lambda\in\R$, we have
\bel{l1c}(-\phi''+V\phi-\lambda\phi=0\textrm{ in }(0,\ell))\quad\Leftrightarrow\quad (-\phi_\star''+V_\star\phi_\star-\ell^2\lambda\phi_\star=0\textrm{ in }(0,1)),\ee
with
$\phi_\star(x)=\phi(\ell x)$ and $ V_\star(x)=\ell^2V(\ell x)$ for $ x\in(0,1).$
Thus, for $j=1,2$, $(\ell^2\lambda_{j,k})_{ k\in\mathbb N}$ are the eigenvalues of the operator $-\partial_x^2+V_{j,\star}$ acting on $L^2(0,1)$ with its domain $H^2(0,1)\cap H^1_0(0,1)$, with
$ V_{j,\star}(x)=\ell^2V_j(\ell x)$ for $ x\in(0,1).$
Combining this with \cite[Theorem 4, p. 35]{PT} gives
$$\ell^2\lambda_{j,k}=(k\pi)^2+\int_0^1V_{j,\star}(x)\d x+\underset{k\to+\infty}{ o}(1)$$
which directly implies \eqref{l1a}. Next we denote by $m_1(\mu)$ the integer
$m_1(\mu):=\sharp \{k\in\mathbb N:\ \ell^2\lambda_{1,k}\leq \mu\ \textrm{and }\int_0^1 \tilde{f}_1(\ell x)\phi_{ 1,k}(\ell x)\d x=0\}.$
By combining \eqref{l1c} with \cite[Lemma 3.3]{Feizmohammadi2025}, we deduce
$$\limsup_{\mu\to+\infty} \frac{m_1(\mu)}{\sqrt{\mu}}\leq \frac{r}{\ell\pi}.$$
Furthermore, we have
$$\int_0^1\tilde{f}_1(\ell x)\phi_{ 1,k}(\ell x)\d x=0\quad \Leftrightarrow\quad \langle \tilde{f}_1, \phi_{ 1,k}\rangle=0,$$
which implies  $m(\mu)=m_1(\ell^{2}\mu)$. Thus we find
$$\limsup_{\mu\to+\infty} \frac{m(\mu)}{\ell\sqrt{\mu}}=\limsup_{\mu\to+\infty}\frac{m_1(\ell^{2}\mu)}{\sqrt{\ell^{2}\mu}}\leq \frac{r}{\ell\pi},$$
and deduce \eqref{l1b}.
\end{proof}

Next we consider an extension of \cite[Theorem A.3]{GS}. To this end, we fix $\psi_1$ a bijection of $\mathbb N$ and
 $$N_j(\mu)=\sharp \{k\in\mathbb N:\ \lambda_{j,k}\leq \mu\},\quad j=1,2,\quad S(\mu)=\sharp \{k\in\mathbb N:\ \lambda_{1,k}\leq \mu,\ \lambda_{1,k}=\lambda_{2,\psi_1(k)}\},\quad \mu>0.$$

\begin{proposition}\label{p4} Let $\alpha\in(0,1)$ and assume that there exists $\mu_0>0$ such that
\bel{p4a}S(\mu) \geq (1-\alpha)\max(N_1(\mu),N_2(\mu))+\frac{\alpha}{2},\quad \mu>\mu_0,\ee
and also that
\bel{p4aa} \mathrm{supp} (V_1-V_2) \subset  \left[0, \frac{(1-\alpha)\ell}{2} \right].\ee
Then we have
\bel{p4b}V_1(x)=V_2(x),\quad x\in\left(0,\ell\right).\ee
\end{proposition}
\begin{proof}
In view of \eqref{p4a},  there exists $\mu_1>0$ such that the condition
$$S(\ell^{-2}\mu)\geq (1-\alpha)\max(N_1(\ell^{-2}\mu),N_2(\ell^{-2}\mu))+\frac{\alpha}{2},\quad \ell^{-2}\mu>\mu_1$$
is fulfilled.
By combining this with \eqref{l1c}, \eqref{p4aa} and \cite[Theorem A.3]{GS}, we deduce
$$\ell^2V_1(\ell x)=\ell^2V_2(\ell x),\quad x\in\left(0,\frac{(1-\alpha)}{2}\right),$$
thus proving the claim after a scaling.\end{proof}

Finally, consider an inverse spectral result involving some extra  knowledge of eigenfunctions of the operator $A_j$, $j=1,2$. Namely, for all $ \mu>0$ and $a\in(0,\ell)$, let
$$E(a,\mu)=\sharp \{k\in\mathbb N:\ \lambda_{1,k}\leq \mu,\ \lambda_{1,k}=\lambda_{2,\psi_1(k)},\ \phi_{1,k}(a)\phi_{2,\psi_1(k)}'(a)=\phi_{2,\psi_1(k)}(a)\phi_{1,k}'(a)\}.$$
By \cite[Theorem 1]{GW}, we have the following result.
\begin{proposition}\label{p5} Let $a_1,a_2\in(0,\ell)$, with $a_1<a_2$, $a_1\leq \frac{\ell}{2}$ and $a_1+a_2\geq\ell$. Let $V_j$, $j=1,2$, be $C^{2k}$ in a neighborhood of $a_1$ and $a_2$, and
also there exists $\mu_0>0$ such that
\begin{align}
\label{p5a}
S(\mu) &\geq 2\frac{a_1}{\ell}\max(N_1(\mu),N_2(\mu))+\left(\frac{1}{2}-\frac{a_1}{\ell}\right)-(k+1),\quad \mu>\mu_0,\\
\label{p5b}E(a_2,\mu)& \geq 2\left(1-\frac{a_2}{\ell}\right)\max(N_1(\mu),N_2(\mu))+\left(\frac{a_2}{\ell}-\frac{1}{2}\right)-(k+1),\quad \mu>\mu_0,
\end{align}
and also that
\bel{p5c} V_1(x)=V_2(x),\quad x\in[a_1,a_2].\ee
Then we have
\bel{p5d}V_1(x)=V_2(x),\quad x\in\left(0,\ell\right).\ee
\end{proposition}
\begin{proof}
By \eqref{p5a}-\eqref{p5b}, we deduce that there exists $\mu_1>0$ such that the conditions
\begin{align*}
S(\ell^{-2}\mu) &\geq \frac{2a_1}{\ell}\max(N_1(\ell^{-2}\mu),N_2(\ell^{-2}\mu))+\left(\frac{1}{2}-\frac{a_1}{\ell}\right)-(k+1),\quad \mu>\mu_1,\\
E(a_2,\ell^{-2}\mu) &\geq 2\left(1-\frac{a_2}{\ell}\right)\max(N_1(\ell^{-2}\mu),N_2(\ell^{-2}\mu))+\left(\frac{a_2}{\ell}-\frac{1}{2}\right)-(k+1),\quad \mu>\mu_1
\end{align*}
are fulfilled.
Combining this with \eqref{l1c}, \eqref{p5c},  and \cite[Theorem 1]{GW} yields
$$\ell^2V_1(\ell x)=\ell^2V_2(\ell x),\quad x\in(0,1),$$
which implies \eqref{p5d} after a scaling.\end{proof}

\bibliographystyle{alpha}
\bibliography{frac}

\newcommand{\etalchar}[1]{$^{#1}$}
\begin{thebibliography}{SQW{\etalchar{+}}20}

\bibitem[AG92]{1}
E~Eric Adams and Lynn~W. Gelhar.
\newblock Field study of dispersion in a heterogeneous aquifer: 2. spatial
  moments analysis.
\newblock {\em Water Res. Research}, 28(12):3293--3307, 1992.

\bibitem[Bru94]{Brusseau:1994}
Mark~L. Brusseau.
\newblock Transport of reactive contaminants in heterogeneous porous media.
\newblock {\em Rev. Geophys.}, 32(3):221--336, 1994.

\bibitem[CJLZ23]{CJLZ}
Siyu Cen, Bangti Jin, Yikan Liu, and Zhi Zhou.
\newblock Recovery of multiple parameters in subdiffusion from one lateral
  boundary measurement.
\newblock {\em Inverse Problems}, 39(10):104001, 31, 2023.

\bibitem[CNYY09]{CJYT}
Jin Cheng, Junichi Nakagawa, Masahiro Yamamoto, and Tomohiro Yamazaki.
\newblock Uniqueness in an inverse problem for a one-dimensional fractional
  diffusion equation.
\newblock {\em Inverse Problems}, 25(11):115002, 16, 2009.

\bibitem[Fei25]{Feizmohammadi2025}
Ali Feizmohammadi.
\newblock Reconstruction of {1-D} evolution equations and their initial data
  from one passive measurement.
\newblock {\em SIAM J. Math. Anal.}, in press, 2025.

\bibitem[GS00]{GS}
Fritz Gesztesy and Barry Simon.
\newblock Inverse spectral analysis with partial information on the potential.
  {II}. {T}he case of discrete spectrum.
\newblock {\em Trans. Amer. Math. Soc.}, 352(6):2765--2787, 2000.

\bibitem[GW15]{GW}
Yongxia Guo and Guangsheng Wei.
\newblock Inverse {S}turm-{L}iouville problems with the potential known on an
  interior subinterval.
\newblock {\em Appl. Anal.}, 94(5):1025--1031, 2015.

\bibitem[HH98]{18}
Yuko Hatano and Naomichi Hatano.
\newblock Dispersive transport of ions in column experiments: An explanation of
  long-tailed profiles.
\newblock {\em Water Res. Research}, 34(5):1027--1033, 1998.

\bibitem[HJK25]{HongJinKian2024}
Jiho Hong, Bangti Jin, and Yavar Kian.
\newblock Identification of a spatially-dependent variable order in
  one-dimensional subdiffusion.
\newblock {\em SIAM J. Math. Anal.}, 57(2):1315--1341, 2025.

\bibitem[HLYZ20]{HLYZ}
Tapio Helin, Matti Lassas, Lauri Ylinen, and Zhidong Zhang.
\newblock Inverse problems for heat equation and space-time fractional
  diffusion equation with one measurement.
\newblock {\em J. Differ. Equations}, 269(9):7498--7528, 2020.

\bibitem[Jin21]{Jin:2021book}
Bangti Jin.
\newblock {\em {Fractional Differential Equations --- An Approach via
  Fractional Derivatives}}, volume 206 of {\em Applied Mathematical Sciences}.
\newblock Springer, Cham, 2021.

\bibitem[JK18]{JK3}
Jaan Janno and Nataliia Kinash.
\newblock Reconstruction of an order of derivative and a source term in a
  fractional diffusion equation from final measurements.
\newblock {\em Inverse Problems}, 34(2):025007, 19, 2018.

\bibitem[JK22]{JK1}
Bangti Jin and Yavar Kian.
\newblock Recovery of the order of derivation for fractional diffusion
  equations in an unknown medium.
\newblock {\em SIAM J. Appl. Math.}, 82(3):1045--1067, 2022.

\bibitem[JLLY17]{JLLY}
Daijun Jiang, Zhiyuan Li, Yikan Liu, and Masahiro Yamamoto.
\newblock Weak unique continuation property and a related inverse source
  problem for time-fractional diffusion-advection equations.
\newblock {\em Inverse Problems}, 33(5):055013, 22, 2017.

\bibitem[JZ21]{JZ}
Bangti Jin and Zhi Zhou.
\newblock Recovering the potential and order in one-dimensional time-fractional
  diffusion with unknown initial condition and source.
\newblock {\em Inverse Problems}, 37(10):105009, 28 pp., 2021.

\bibitem[KJ19]{JK}
N~Kinash and Jaan Janno.
\newblock An inverse problem for a generalized fractional derivative with an
  application in reconstruction of time- and space-dependent sources in
  fractional diffusion and wave equations.
\newblock {\em Mathematics}, 7:1138, 2019.

\bibitem[KLY23]{KLY23}
Yavar Kian, Yikan Liu, and Masahiro Yamamoto.
\newblock Uniqueness of inverse source problems for general evolution
  equations.
\newblock {\em Commun. Contemp. Math.}, 25(6):2250009, 33, 2023.

\bibitem[Kou08]{57}
Samuel~C Kou.
\newblock Stochastic modeling in nanoscale biophysics: subdiffusion within
  proteins.
\newblock {\em Ann. Appl. Stat.}, 2(2):501--535, 2008.

\bibitem[KR23]{KaltenbacherRundell:2023book}
Barbara Kaltenbacher and William Rundell.
\newblock {\em Inverse problems for fractional partial differential equations}.
\newblock AMS, Providence, RI, 2023.

\bibitem[KRY20]{KRY}
Adam Kubica, Katarzyna Ryszewska, and Masahiro Yamamoto.
\newblock {\em {Time-Fractional Differential Equations---A Theoretical
  Introduction}}.
\newblock Springer, Singapore, 2020.

\bibitem[KSXY22]{KSXY}
Yavar Kian, \'Eric Soccorsi, Qi~Xue, and Masahiro Yamamoto.
\newblock Identification of time-varying source term in time-fractional
  evolution equations.
\newblock {\em Commun. Math. Sci.}, 20(1):53--84, 2022.

\bibitem[KY17]{KY1}
Yavar Kian and Masahiro Yamamoto.
\newblock On existence and uniqueness of solutions for semilinear fractional
  wave equations.
\newblock {\em Fract. Calc. Appl. Anal.}, 20(1):117--138, 2017.

\bibitem[Lev44]{Levenberg:1944}
Kenneth Levenberg.
\newblock A method for the solution of certain non-linear problems in least
  squares.
\newblock {\em Quart. Appl. Math.}, 2:164--168, 1944.

\bibitem[LIY16]{LiImanuvilovYamamoto:2016}
Zhiyuan Li, Oleg~Yu. Imanuvilov, and Masahiro Yamamoto.
\newblock Uniqueness in inverse boundary value problems for fractional
  diffusion equations.
\newblock {\em Inverse Problems}, 32(1):015004, 16, 2016.

\bibitem[Mar63]{Marquardt:1963}
Donald~W. Marquardt.
\newblock An algorithm for least-squares estimation of nonlinear parameters.
\newblock {\em J. Soc. Indust. Appl. Math.}, 11:431--441, 1963.

\bibitem[MJCB14]{72}
Ralf Metzler, Jae-Hyung Jeon, Andrey~G Cherstvy, and Eli Barkai.
\newblock Anomalous diffusion models and their properties: non-stationarity,
  non-ergodicity, and ageing at the centenary of single particle tracking.
\newblock {\em Phys. Chem. Chem. Phys.}, 16(44):24128--24164, 2014.

\bibitem[MY13]{MillerYamamoto:2013}
Luc Miller and Masahiro Yamamoto.
\newblock Coefficient inverse problem for a fractional diffusion equation.
\newblock {\em Inverse Problems}, 29(7):075013, 8, 2013.

\bibitem[Nig86]{76}
Raoul~R. Nigmatullin.
\newblock The realization of the generalized transfer equation in a medium with
  fractal geometry.
\newblock {\em Phys. Stat. Sol. B}, 133:425--430, 1986.

\bibitem[NSY10]{NSY}
Junichi Nakagawa, Kenichi Sakamoto, and Masahiro Yamamoto.
\newblock Overview to mathematical analysis for fractional diffusion equations
  -- new mathematical aspects motivated by industrial collaboration.
\newblock {\em J. Math-for-Ind.}, 2(A):99--108, 2010.

\bibitem[Pod99]{P}
Igor Podlubny.
\newblock {\em Fractional {D}ifferential {E}quations}.
\newblock Academic Press, Inc., San Diego, CA, 1999.

\bibitem[PT87]{PT}
J\"urgen P\"oschel and Eugene Trubowitz.
\newblock {\em Inverse {S}pectral {T}heory}.
\newblock Academic Press, Inc., Boston, MA, 1987.

\bibitem[SQW{\etalchar{+}}20]{Sun:2024water}
Liwei Sun, Han Qiu, Chuanhao Wu, Jie Niu, and Bill~X. Hu.
\newblock A review of applications of fractional
  advection–dispersionequations for anomalous solute transport in surface
  andsubsurface water.
\newblock {\em WIREs Water}, 7(4):e1448, 21 pp., 2020.

\bibitem[Yam23]{Y}
Masahiro Yamamoto.
\newblock Uniqueness for inverse source problems for fractional diffusion-wave
  equations by data during not acting time.
\newblock {\em Inverse Problems}, 39(2):024004, 20, 2023.

\end{thebibliography}

\end{document}